\documentclass{amsart}
\usepackage{amsmath,amssymb}
\usepackage[dvips]{graphicx}
\usepackage{array}

\newtheorem{thm}{Theorem}
\newtheorem{prop}[thm]{Proposition}
\newtheorem{lem}[thm]{Lemma}

\newtheorem{defi}[thm]{Definition}

\newtheorem*{notat}{Notation}
\newtheorem*{spthm}{Main Theorem}

\newcommand{\Cg}{{\mathbb C}}

\newcommand{\Ng}{{\mathbb N}}

\newcommand{\Rg}{{\mathbb R}}

\newcommand{\Zg}{{\mathbb Z}}

\newcommand{\Pf}{\noindent{\bf Proof.}}
\newcommand{\Cr}{{\mathcal C}}
\newcommand{\Co}{{\widehat{\mathcal X}}}
\newcommand{\Cac}{{\mathcal X}}

\newcommand{\Or}{{\mathcal O}}

\newcommand{\Mar}{{\mathcal M}}
\newcommand{\Nar}{{\mathcal N}}

\begin{document}
\title[Minimizing Parabolic orbits]{Globally minimizing parabolic motions in the \\
Newtonian $N$-body Problem}

\author{E. Maderna \and A. Venturelli}
\address{E. Maderna:  Instituto de Matem\'atica y Estad\'istica ``Prof. Rafael 
Laguardia'',
Universidad de la Rep\'ublica,  
Herrera y Reissig 565, 12000 Montevideo (UY)}
\address{A. Venturelli: Laboratoire d'Analyse non lin\'eaire et G\'eometrie
Universit\'e d'Avignon et des pays de Vaucluse
33, Rue Louis Pasteur - 84000, Avignon (FR)}      

\email {Andrea.Venturelli@univ-avignon.fr, eze@fing.edu.uy }
\maketitle
\date{}

%
%
\begin{abstract}
We consider the $N$-body problem in $\Rg^d$ with the newtonian  
potential $1/r$. We prove that for every initial configuration $x_i$ and for every minimizing 
normalized central configuration $x_0$, there exists a collision-free parabolic solution starting from 
$x_i$ and asymptotic to $x_0$. This solution is a minimizer 
in every time interval.  
The proof exploits the variational structure of the problem, and it consists in finding a convergent 
subsequence in a family of minimizing trajectories. The hardest part is to show that this 
solution is parabolic and asymptotic to $x_0$. 
\end{abstract}
\vspace{2mm}
\section{\bf Introduction}
\vspace{3mm}
%
%
%
%
\vspace{2mm}
In this paper we consider $N$ positive masses in an euclidean space $\Rg^d$, submitted to a gravitational 
interaction. We find some interesting solutions with a given asymptotic behaviour. 
The equation of motion of the $N$-body problem is written  
\begin{equation} \label{Newt}
\ddot{\vec{r}}_i=-\sum\limits_{j=1,...,N,\ j\neq i} \frac{m_j 
(\vec{r}_i -\vec{r}_j)}{|\vec{r}_i -\vec{r}_j|^3}.
\end{equation}
where $m_i$ is the mass and $\vec{r}_i\in\Rg^d$ the position of the $i$-th body.
Since these equations are invariant by translation, we can assume that the center of mass is 
at the origin. \\ \noindent 
These equations are Euler-Lagrange equations of the Lagrangian action functional (we will define 
it precisely in the next section), therefore solutions of (\ref{Newt}) are critical points of the action
in a set of paths with fixed ends. The simplest kind of critical points are minima, so it is natural to search
for minimizers of the lagrangian action joining two given configurations in a fixed time.  
The potential of the N-body problem being singular at collision configurations, a main difficult 
involved in this approach is to show that minimizers are collision-free. The following theorem, essentially 
due to C. Marchal, is a major advanced in this subject. 
\begin{thm}
Given two N-body configurations $x_i=(\vec{r}_1,...,\vec{r}_N)\in (\Rg^d)^N$, 
$x_f=(\vec{s}_1,...,\vec{s}_N)\in (\Rg^d)^N$
and a time $T>0$, an action minimizing path joining $x_i$ to $x_f$ in time $T$ is collision-free for $t\in (0,T)$. 
\end{thm}
\noindent
See [Ma1], [Ma2] and [Ch2] for a claim and a proof of this theorem for $d=2$ and $d=3$. See [Fe-Te] for a proof 
in any dimension. This theorem, together with the lower semicontinuity of the action (see Section \ref{variat-setting}), 
implies in particular that there always exists a collision-free minimizing solution joining two given collision-free N-body 
configuration in a given time.  \\ \noindent
A natural extension of Marchal's theorem is to search solutions defined on an infinite interval $[0,+\infty)$, starting
from a given configuration at $t=0$ and having a given asymptotic behaviour for $t\rightarrow +\infty$. 
The classification of all possible asymptotic behaviour of solutions in the N-body problem has been investigated since the 
beginning of the last century. The main results in this direction are due to J. Chazy. In [Cha1] it is 
shown that there are only seven possible final evolutions in the three-body problem. 
Among these seven possibilities there are the so-called {\sl parabolic motions}. 
A solution $t\mapsto (\vec{r}_1,...,\vec{r}_N)(t)$ of the N-body
problem is said to be parabolic if the velocity of every body tends to zero as $t\rightarrow +\infty$.     
We introduce the functions 
\begin{equation} \label{potentiel&moment-d-inertie}
I(x)=\sum\limits_{i=1}^N m_i |\vec{r}_i|^2,\qquad U(x)=\sum\limits_{1\le i<j\le N} \frac{m_i m_j}{|\vec{r}_i-\vec{r}_j|},
\qquad x=(\vec{r}_1,...,\vec{r}_N),
\end{equation}
respectively equal to the moment of inertia with respect to the center of mass and to the Newtonian potential. 
\begin{notat}
Given a configuration $x$, we denote by $\tilde{x}=I(x)^{-1/2} x$ the associated normalized configuration.
\end{notat}
\noindent
It is well known (see for istance [Hu-Sa] and [Ch1]) that if $t\mapsto x(t)$ is a parabolic solution, 
the normalized trajectory $\tilde{x}(t)$ is asymptotic to the set of {\sl central configurations} 
(i.e. critical points of $\tilde{U}=I^{1/2}U$).
Given a central configuration $x_0$ with $I(x_0)=1$, we say that a parabolic solution $t\mapsto x(t)$ 
is asymptotic to $x_0$ if $\tilde{x}(t)\rightarrow x_0$ as $t\rightarrow +\infty$. 
A central configuration $x_0$ is said to be minimizing if it is an absolute minimum of $\tilde{U}$. 
We can now state the main result of this paper.
\begin{spthm}
Given any initial configuration $x_i$ and any minimizing normalized central configuration $x_0$, 
there exists a parabolic solution $\gamma : [0,+\infty)\rightarrow (\Rg^d)^N$ 
starting from $x_i$ at $t=0$ and asymptotic to $x_0$ for
$t\rightarrow +\infty$. This solution is a minimizer of the lagrangian action with fixed ends 
in every compact interval contained in $[0,+\infty)$ and it is collision-free for $t>0$.  
\end{spthm}   
\noindent
We do not require any hypothesis of nondegeneracy of the central configuration $x_0$. 
\\ \noindent
The parabolic solution $\gamma$ is constructed as limit of a sequence $\gamma_n :[0,t_n]\rightarrow (\Rg^d)^N$
of minimizers connecting $x_i$ with a configuration homothetic to $x_0$ in time $t_n$, and $t_n\rightarrow +\infty$.
In Section \ref{construct-solut} we construct the sequence $\gamma_n$ and we prove that it is uniformly convergent on every 
compact subset of $\Rg$.  
In Sections \ref{parabol-solut} and \ref{Lambert} we show that $\gamma$ is parabolic and asymptotic to $x_0$. 
The proof of this last
property is achieved by comparing the action of the N-body problem with the action of a Kepler problem, and using
Lambert's Theorem to estimates the action. 
In the Appendix we state and prove some technical estimates concerning the Kepler problem on the line
that we need to construct $\gamma$ and to prove its parabolicity.  
The authors believe that these minimizing parabolic solutions are in fact calibrated 
curves of some weak KAM solutions of the N-body problems, whose existence has been proved in [Mad] by one of 
the authors. 
Our Main Theorem has a natural interpretation in terms of McGehee vector field and collision manifold.
Indeed, in [Ch1], [McG] and [Mo] it is shown that if $x_0$ is a central configuration with $I(x_0)=1$,  
the state $(x_0, v_0 x_0)$ with $v_0=(2 U(x_0)^{1/2}$ is a critical point 
of the McGehee vector field in the collision manifold, and 
its stable set corresponds to parabolic solutions asymptotic to $x_0$ as $t\rightarrow +\infty$. Thus, we can 
formulate the Main Theorem by saying that the stable set of $(x_0, v_0 x_0)$ (for the McGehee vector field) projects on the 
whole configuration space, provided $x_0$ is a minimizing central configuration. \\ \noindent
We think that variational methods could be used to study some important features on the global dynamics
of N-body problem. In particular, it should be interesting to study hyperbolic solutions using 
variational methods. We recall that a solution $\gamma :[0,+\infty)\rightarrow (\Rg^d)^N$ is said to be hyperbolic if
there exists a (collision-free) configuration $x_0$ such that 
\begin{equation} \label{solut-hyperb}
\gamma(t)=x_0 t+o(t), \qquad t\rightarrow +\infty
\end{equation}
A hyperbolic solution has necessarily positive energy, and replacing $x_0$ by a normalized configuration,
(\ref{solut-hyperb}) is equivalent to $\gamma(t)=\sqrt{2h} x_0 t+o(t)$ as $t\rightarrow +\infty$ (see [Cha1]), 
where $h$ is the energy of the solution.  
In this case we will say that $\gamma(t)$ is hyperbolic for $t\rightarrow +\infty$ and asymptotic to $x_0$.   
Since there is no constraint to the limit configuration $x_0$ of a hyperbolic solution (see again [Cha1]), 
it is natural to ask the following two questions. The second one has been asked by R. Montgomery. 
\\ \par
\noindent
{\bf Question 1.} Given an initial configuration $x_i$ and a normalized non-collision configuration $x_0$, 
does there exist a hyperbolic motion  
starting from $x_i$ at $t=0$ and asymptotic to $x_0$ for $t\rightarrow +\infty$ ?
\vspace{1mm} \\ \noindent
{\bf Question 2.} For which couple of normalized non-collision configurations $x_0$ and $x_0^\prime$ does there exist a 
solution that is hyperbolic both for $t\rightarrow +\infty$ and for 
$t\rightarrow -\infty$ and is asymptotic to $x_0$ for $t\rightarrow +\infty$ and to $x_0^\prime$ for $t\rightarrow -\infty$ ? 
\vspace{1mm} \\ \noindent
We hope that it will be possible to answer these questions using variational methods similar to those 
developed in this paper.  
%
%
%
%
\vspace{2mm}
\section{\bf Variational setting} \label{variat-setting}
\vspace{3mm} \noindent
Since equations (\ref{Newt}) are invariant by translation, we fix the origin of our 
inertial frame at the center of mass of the system.  
We define the {\sl configuration space} of the 
system as 
$$
\Cac=\left\{ x=(\vec{r}_1,...,\vec{r}_N)\in (\Rg^d)^N,\quad 
\sum\limits_{i=1}^N m_i\vec{r}_i=0 \right\},
$$
and we endow $\Cac$ with the {\sl mass scalar product} : 
$$
\begin{array}{rl}
x\cdot y&=\sum\limits_{i=1}^N m_i <\vec{r}_i,\vec{s}_i> \\
x&=(\vec{r}_1,...,\vec{r}_N)\in\Cac,\quad y=(\vec{s}_1,...,\vec{s}_N)\in\Cac,
\end{array}
$$
where $<\ ,\ >$ is the usual euclidean product in $\Rg^d$. 
We denote by $\|\ \|$ the euclidean norm on $\Cac$ associated to the mass scalar product.
A configuration $x=(\vec{r}_1,...,\vec{r}_N)\in \Cac$ is said to be a 
{\sl collision configuration} if $\vec{r}_i=\vec{r}_j$ for some $i\neq j$. 
We denote by $Coll$ the set of collision configurations and  
by $\Co=\Cac\setminus Coll$ the set of collisions-free configurations. 
Equations (\ref{Newt}) can be written in a more compact form as a 
second order differential equation on $\Co$ 
\begin{equation} \label{Newt-2}
\ddot{x}=\nabla U(x),   
\end{equation}
where $U$ is the newtonian potential already defined in (\ref{potentiel&moment-d-inertie}), 
the gradient is calculated with 
respect to the mass scalar product. Since $\Co$ is an open subset of $\Cac$, 
the tangent space of $\Co$ is identified with $\Co\times \Cac$. The following 
functions defined on $\Co\times \Cac$
$$
K=y\cdot y,\quad L=\frac{K}{2}+U,\quad H=\frac{K}{2}-U,
$$
are respectively equal to twice the 
kinetic energy, to the lagrangian and to the energy first integral.  
\\ \noindent
Given an absolutely continuous path $\gamma :[a,b]\rightarrow \Cac$, 
we define its {\sl Lagrange action} by :
$$
A_L(\gamma)=\int_a^b L(\gamma(t),\dot{\gamma}(t))dt,
$$
where $L$ is naturally extended to a function defined over $\Cac\times \Cac$ by $L(x,y)=+\infty$ 
if $x\in Coll$.
It is well known that collision-free extremals 
of $A_L$ are solutions of equations (\ref{Newt-2}). 
\begin{defi}
We say that an absolutely continuous path $\gamma :[a,b]\rightarrow \Cac$ is a minimizer if 
$A_L(\sigma)\ge A_L(\gamma)$ for every absolutely continuous path $\sigma :[a,b]\rightarrow \Cac$
having the same extremities. 
If $I\subset \Rg$ is any interval, we say that 
$\gamma : I\rightarrow \Cac$ is minimizing if for every compact interval $[a,b]$ contained in $I$,
the path $\gamma \left|_{[a,b]}\right.$ is a minimizer.   
\end{defi} \noindent
Given a positive real number $T$ and two configurations $x_i$ and $x_f$, let 
$\Sigma(x_i,x_f;T)$ be the set of absolutely continuous paths defined in the interval $[0,T]$ and joining 
$x_i$ to $x_f$ in time $T$. The following proposition is well known.
\begin{prop}
For every $x_i,x_f \in \Cac$ and for every $T>0$ there exists a minimizer 
$\gamma :[0,T]\rightarrow \Cac$
joining $x_i$ to $x_f$. 
\end{prop}
\noindent
In [Ve] and [Fe-Te] one can find a proof of this proposition when the functional $A_L$ is defined 
over $H^1$ paths (i.e. absolutely continuous paths with derivative in $L^2$.)
joining $x_i$ to $x_f$. An absolutely continuous path having a finite action is necessarily in $H^1$,
therefore minimizers among $H^1$ paths are also minimizers among absolutely continuous paths.
\vspace{1mm}\\ \noindent
The proposition above do not ensure that $\gamma$ is collision-free, but by the already cited Marchal's theorem, if 
$d\ge 2$, minimizers are collision-free for $t\in (0,T)$.  
%
%
%
%
\vspace{2mm}
\section{\bf Construction of the solution} \label{construct-solut}
\vspace{3mm} \noindent
In this section we construct the solution $\gamma :[0,+\infty)\rightarrow \Cac$ of the main theorem as 
limit of minimizers. We will show in Sections \ref{parabol-solut} and \ref{Lambert} that $\gamma$ is parabolic and asymptotic 
to $x_0$. \\ \noindent
Before stating the main result of this section, we recall a classical result concerning parabolic solutions 
(see [Ch1] or [Hu-Sa]) for a proof). 
\begin{prop} \label{parabolic_sundman}
If $\gamma :[0,+\infty)\rightarrow \Cac$ is a parabolic solution of the N-body problem, the energy of $\gamma$ is necessarily
zero, moreover we have 
$$
I(t) =\alpha^2 t^{\frac{4}{3}}+o(t^{\frac{4}{3}}), \qquad \nabla \tilde{U}\left(\tilde{\gamma}(t)\right) 
\rightarrow 0, \qquad \tilde{U}\left(\tilde{\gamma}(t)\right)\rightarrow U_0
$$
as $t\rightarrow +\infty$, where 
\begin{equation} \label{defit-alpha}
\alpha=\left(9 U_0/2\right)^{1/3}
\end{equation}
In particular, the $\omega$-limit of $\tilde{\gamma}(t)$ 
is contained in the set of normalized central configuration. 
\end{prop} 
\noindent
Since there are always infinitely many normalized central configurations for a given critical level of $\tilde{U}$, 
(the orthogonal group acts on $\Co$ leaving invariant $\tilde{U}$), we cannot say 
{\it a priori} that the $\omega$-limit of $\tilde{\gamma}(t)$ is a given configuration.
If $\gamma(t)$ is a parabolic solution asymptotic to normalized central configuration $x_0$ 
(i.e. $\tilde{\gamma}(t)$ converges to $x_0$), by Proposition \ref{parabolic_sundman} we have the asymptotic estimates
\begin{equation} \label{gamma-I_0-t43}
\gamma(t)= \alpha x_0 t^\frac{2}{3}+o(t^\frac{2}{3}),\quad \text{as}\quad t\rightarrow +\infty
\end{equation}
The following Lemma is a converse of Proposition \ref{parabolic_sundman}. 
\begin{lem} \label{lem-parab-asympt-x0}
Let $x_0$ be a normalized central configuration, $U_0=\tilde{U}(x_0)$ and $\alpha$ the constant 
defined in (\ref{defit-alpha}).
A solution $\gamma : [0,+\infty)\rightarrow \Cac$ satisfying the asymptotic estimates (\ref{gamma-I_0-t43}) is parabolic and 
asymptotic to $x_0$.
\end{lem}
\Pf\  We just need to prove that $\gamma$ is parabolic. Replacing (\ref{gamma-I_0-t43}) in the equation of motion we find
$\ddot{\gamma}(t)=\Or(t^{-\frac{4}{3}})$, as $t\rightarrow +\infty$. Therefore, the velocity $\dot{\gamma}(t)$ has 
a limit for $t\rightarrow +\infty$ that we denote $\dot{\gamma}_\infty$. Moreover we have 
$$
\dot{\gamma}(t)=\dot{\gamma}_\infty+\Or(t^{-\frac{1}{3}}),\qquad t\rightarrow +\infty.
$$
Integrating this expression we find 
$$
\gamma(t)=\dot{\gamma}_\infty t +\Or(t^\frac{2}{3}),\qquad t\rightarrow +\infty,
$$
thus $\dot{\gamma}_\infty=0$ and 
$\gamma(t)$ is parabolic. 
\vspace{1mm} \\ \noindent    
By the way, if $x_0$ is a normalized central configuration, the path 
\begin{equation} \label{pauline}
\gamma_0 : [0,+\infty)\rightarrow \Cac,\qquad \gamma_0(t)=\alpha x_0 t^\frac{2}{3},
\end{equation}
is a solution of the $N$-body problem. $\gamma_0$ is called {\sl homothetic-parabolic solution asymptotic to} $x_0$. 
\\ \noindent
We state now the main result of this section. We recall that $x_i$ is the initial configuration of the Main Theorem, 
$x_0$ is a normalized minimizing central configuration, $U_0$ and $\alpha$ are as before, $\gamma_0(t)$ is given by 
(\ref{pauline}). 
\begin{thm} \label{thm-exist-solut}
There exists a minimizing solution $\gamma :[0,+\infty)\rightarrow \Cac$ starting from $x_i$, 
a sequence of positive numbers $t_n\rightarrow +\infty$ and a sequence of minimizers
$\gamma_n\in \Sigma(x_i,\gamma_0(t_n);t_n)$ such that $\gamma_n$ converges uniformly to $\gamma$ 
on every compact interval contained in $[0,+\infty)$. Moreover $\gamma(t)$ is collision-free for $t>0$.
\end{thm}
\noindent
We prove this Theorem in several steps. At Proposition \ref{lacapelle-marival} we show that
if $T$ and $t/T$ are sufficiently great, for every minimizer $\overline{\gamma}\in\Sigma(x_i,\gamma_0(t);t)$ the action
$A_L(\overline{\gamma}\left|_{[0,T]}\right.)$ has a uniform bound (independent of $t$). Successively, using Ascoli's theorem and 
a diagonal trick, we find the sequence $(\gamma_n)_{n=1}^{+\infty}$.   
We start with some preliminary definitions and remarks.
Given two configurations $x$ and $x^\prime$ and a time $T$, we denote by ${\mathcal A}(x,x^\prime;T)$ the action of a 
minimizing path joining $x$ to $x^\prime$ in time $T$ (the same function is denoted $\phi(x,x^\prime,T)$ in [Mad]). 
In a similar way, given two positive real numbers $a$ and $b$ 
and a time $T$, we denote by $S(a,b;T)$ the action (for the one dimensional keplerian problem with 
lagrangian $\frac{\dot{r}^2}{2}+\frac{U_0}{r}$) of a minimizing path joining $a$ to $b$ in time $T$. 
\\ \noindent
By the homogeneity of $U$, if $\varpi :[0,T]\rightarrow \Cac$ is a solution of (\ref{Newt-2}) and $\lambda>0$, the path 
$$
\varpi^\lambda :[0,\lambda T]\rightarrow \Cac,\qquad \varpi^\lambda(t)=\lambda^{\frac{2}{3}}\varpi(t/\lambda)
$$
is still a solution of (\ref{Newt-2}). Moreover, if $\varpi$ is a minimizer, $\varpi^\lambda$ is still 
a minimizer. A similar property holds for solutions and minimizers of a one dimensional 
Kepler problem. Therefore we have 
$$
{\mathcal A}(\lambda^\frac{2}{3}x, \lambda^\frac{2}{3}x^\prime;\lambda T)=\lambda^\frac{1}{3} {\mathcal A}(x,x^\prime;T),
\qquad S(\lambda^\frac{2}{3}a,\lambda^\frac{2}{3}b;\lambda T)=\lambda^\frac{1}{3} S(a,b;T). 
$$ 
\begin{lem} \label{nonna-isabella}
We have
$$
{\mathcal A}(x,x^\prime;T)\ge S(\|x\|,\|x^\prime\|;T)
$$
with equality if and only if $x$ and $x^\prime$ are on the half-line starting from zero 
generated by $\hat{x}$, where $\hat{x}$ 
is a normalized minimizing configuration (i.e. $\|\hat{x}\|=1$ and $\tilde{U}(\hat{x})=U_0$). 
\end{lem}
\Pf\ 
Let $\varpi :[0,T]\rightarrow \Cac$ be a minimizer joining $x$ to $x^\prime$ in time $T$ and let 
$r(s)=\|\varpi(s)\|$. By Sundman inequality we have
$$
\|\dot{\varpi}(s)\|^2\ge \dot{r}^2(s),
$$
with equality if and only if $\dot{\varpi}(s)$ is parallel to $\varpi(s)$. Since $U_0$ is the minimum 
of $\tilde{U}$ we have also
$$
U(\varpi(s))\ge \frac{U_0}{r(s)},
$$
with equality if and only if $\tilde{U}(\varpi(s))=U_0$. Therefore 
\setlength{\extrarowheight}{0.5cm}
\begin{displaymath}
\begin{array}{rl}
{\mathcal A}(x,x^\prime;T)&=\displaystyle\int_0^T \left( \frac{\|\dot{\varpi}(s)\|^2}{2}+U(\varpi(s))\right) ds \\
&\ge\displaystyle\int_0^T \left( \frac{\dot{r}^2(s)}{2}+\frac{U_0}{r(s)}\right)ds \\
&\ge S(\|x\|,\|x^\prime\|;T)
\end{array}
\end{displaymath}
\setlength{\extrarowheight}{0cm}
with equality if and only if $\varpi(s)=\mu(s)\hat{x}$, where $\hat{x}$ is a minimizing
normalized configuration and $s\mapsto \mu(s)\in\Rg_+$ is a minimizer 
(for the one-dimensional Kepler problem) joining $\|x\|$ to $\|x^\prime\|$ 
in time $T$. This proves the Lemma. 
\vspace{1mm} \\ \noindent
In order to simplify the exposition we introduce the following notation. If $x,x^\prime\in \Cac$ are two 
configurations and $0\le \tau<T<t$ we term
\begin{equation} \label{diff-action}
\Mar(x,x^\prime;\tau,T,t)={\mathcal A}(0,x;T+\tau)+{\mathcal A}(x,x^\prime;t-T)-{\mathcal A}(0,x^\prime;t-\tau).
\end{equation}
In a similar way, if $r,r\prime\in [0,+\infty)$  and $0\le \tau<T<t$ we term
\begin{equation} \label{diff-act-kepl}
\Nar(r,r^\prime;\tau,T,t)=S(0,r;T+\tau)+S(r,r^\prime;t-T)-S(0,r^\prime;t-\tau).
\end{equation}
\begin{lem} \label{lemme-fondamentale}
Let $0<\tau<T<t$ be real numbers. If $\xi\in\Sigma(0,x_i;\tau)$  
and $\overline{\gamma}\in \Sigma(x_i,\gamma_0(t);t)$ we have 
$$
\Mar(\overline{\gamma}(T),\gamma_0(t);\tau,T,t)\le 2 A_L(\xi) \quad \text{and} \quad 
\Nar(\|\overline{\gamma}(T)\|,\alpha t^\frac{2}{3};\tau,T,t)\le 2 A_L(\xi).
$$
\end{lem}
\Pf \
To prove the first inequality, let $\eta\in\Sigma(0,\overline{\gamma}(T);T+\tau)$   
and let $\zeta\in\Sigma(0,\gamma_0(t);t-\tau)$. The path $\zeta$ is nothing but a repametrization
of $\gamma_0\left|_{[0,t]}\right.$.   
Since $\eta$ and $\overline{\gamma}$ are minimizers, we have the triangular inequalities  
$$
\begin{array}{rl}
A_L(\eta)&\le A_L(\xi)+A_L(\overline{\gamma}\left|_{[0,T]}\right.) \\
A_L(\overline{\gamma}) &\le A_L(\xi)+A_L(\zeta),
\end{array}
$$
therefore
$$
A_L(\eta)+A_L(\overline{\gamma}\left|_{[T,t]}\right.)\le A_L(\xi)+A_L(\overline{\gamma}) \le 2 A_L(\xi)+A_L(\zeta).
$$
This gives the first inequality. \\ \noindent
The second inequality is a direct consequence of the first one and of Lemma \ref{nonna-isabella}.
\begin{figure}
\centering
\includegraphics[scale=1.0]{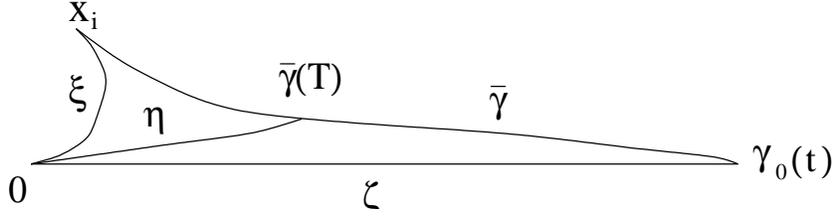}
\caption{The paths $\overline{\gamma}$, $\eta$, $\zeta$ and $\xi$ in the configuration space. }
\label{construct}
\end{figure}
\begin{prop} \label{lacapelle-marival}
There exist three constants $K>0$, $\overline{T}>0$ and $\overline{s}>1$ such that for every
$T\ge \overline{T}$, for every $t\ge \overline{s}T$ and for every $\overline{\gamma}\in  \Sigma(x_i,\gamma_0(t);t)$ 
we have
$$
\|\overline{\gamma}(T)\|\le K T^\frac{2}{3}.
$$
\end{prop}
\Pf\  
Suppose, for the sake of a contradiction, that there exist three sequences 
of positive real numbers $(K_n)_{n=0}^{+\infty}$, $(T_n)_{n=0}^{+\infty}$ and $(t_n)_{n=0}^{+\infty}$ 
satisfying
$$
K_n\rightarrow +\infty,\qquad T_n\rightarrow +\infty,\qquad \frac{t_n}{T_n}\rightarrow +\infty,
$$
and a sequence of minimizers $\overline{\gamma}_n\in \Sigma(x_i,\gamma_0(t_n);t_n)$ such that for every $n\in \Ng$ :
$$
\|\overline{\gamma}_n(T_n)\|\ge K_n T_n^\frac{2}{3}.
$$   
Let $\tau>0$ and $\xi: [0,\tau]\rightarrow \Cac$ be a minimizer connecting $0$ to $x_i$ in time $\tau$.
Without loss of generality we can assume $0<\tau<T_n<t_n$. By homothety invariance 
and by the second inequality of Lemma \ref{lemme-fondamentale} we have
\begin{equation} \label{maria-rosa}
T_n^\frac{1}{3}\Nar\left(\frac{\|\overline{\gamma}_n(T_n)\|}{T_n^\frac{2}{3}},\alpha\left(\frac{t_n}{T_n}\right)^\frac{2}{3}
;\frac{\tau}{T_n},1,\frac{t_n}{T_n}\right)\le 2 A_L(\xi).
\end{equation}
Since $\frac{\|\overline{\gamma}_n(T_n)\|}{T_n^\frac{2}{3}}\rightarrow +\infty$ and $\frac{t_n}{T_n}\rightarrow +\infty$, 
by Proposition \ref{elenuccia} of the Appendix we have 
$$
\Nar\left(\frac{\|\overline{\gamma}_n(T_n)\|}{T_n^\frac{2}{3}},\alpha\left(\frac{t_n}{T_n}\right)^\frac{2}{3}
;\frac{\tau}{T_n},1,\frac{t_n}{T_n}\right)\rightarrow +\infty
$$
as $n\rightarrow +\infty$. This contradicts inequality (\ref{maria-rosa}).
\vspace{1mm}\\ \noindent
We need now an estimates of the minimal action ${\mathcal A}(x,x^\prime;T)$ 
when $\|x\|$ and $\|x^\prime\|$ are less then a given size. 
\begin{prop} \label{cap-de-l-homy}
There exist two positive constants $C_1$ and $C_2$ such that if $R>0$ and $T>0$, if $x\in \Cac$ and 
$x^\prime\in \Cac$ are two configurations satisfying $\|x\|\le R$ and $\|x^\prime\|\le R$, 
we can find an absolutely continuous path 
$\gamma_{x x^\prime} :[0,T]\rightarrow \Cac$ joining $x$ to $x^\prime$ in time $T$ such that the 
following inequality holds 
\begin{equation} \label{ezequiel}
A_L(\gamma_{x x^\prime})\le C_1 \frac{R^2}{T}+C_2 \frac{T}{R}.
\end{equation}
In particular we have 
\begin{equation} \label{budoni}
{\mathcal A}(x,x^\prime;T)\le  C_1 \frac{R^2}{T}+C_2 \frac{T}{R}.
\end{equation}
\end{prop}
 \Pf \
Let $x_0^\prime\in \Co$ be any normalized collision-free configuration. 
We construct an absolutely continuous path $\gamma_x : [0,T/2]\rightarrow \Cac$ joining $x$ to 
$R x_0^\prime$ and verifying
\begin{equation} \label{PMA}
A_L(\gamma_x)\le A_1 \frac{R^2}{T}+A_2 \frac{T}{R}
\end{equation}
where $A_1$ and $A_2$ are two positive constants independent on $R$, $T$ and $x$. An analogous
path $\gamma_{x^\prime}: [0,T/2]\rightarrow \Cac$ joining $R x_0^\prime$ to $x^\prime$ can be 
constructed in exactly the same way. Pasting $\gamma_x$ and $\gamma_{x^\prime}$ together and choosing 
$C_1=2A_1$ and $C_2=2A_2$ we get a path $\gamma_{x x^\prime}$ verifying (\ref{ezequiel}).
Inequality (\ref{budoni}) is an obvious consequence of (\ref{ezequiel}). 
\\ \noindent
Let $x_0^\prime=(\vec{c}_1,...,\vec{c}_N)$. We term $\vec{c}_{ij}=\vec{c}_j -\vec{c}_i$ 
and $c_{ij}=|\vec{c}_{ij}|$ for $1\le i<j\le N$. 
In a similar way, given $x=(\vec{r}_1,...,\vec{r}_N)\in \Cac$ 
with $\|x\|\le R$, we term $\vec{r}_{ij}=\vec{r}_j -\vec{r}_i$ and $r_{ij}=|\vec{r}_{ij}|$. 
Let $\lambda_{ij}$ be the coefficients
$$
\lambda_{ij}=\frac{r_{ij}}{Rc_{ij}+r_{ij}}\in [0,1),
$$ 
and let $h$ be the cardinality of the set $\{\lambda_{ij}\}_{1\le i<j\le N}$.
The inequality $1\le h\le N(N-1)/2$ holds. Let us denote
$$
0\le \mu_1<...<\mu_h<1
$$
the elements of the set $\{\lambda_{ij}\}_{1\le i<j\le N}$ ordered increasingly. 
We define $\mu_0=0$ and $\mu_{h+1}=1$. For every $i=0,...,h$ we term
$$
\tau_i=\frac{T(\mu_{i+1}-\mu_i)^\frac{3}{2}}{2\sum\limits_{k=0}^h (\mu_{k+1}-\mu_k)^\frac{3}{2}}.
$$ 
We observe that $\tau_0\ge 0$ and $\tau_i>0$ if $i\ge 1$, moreover $\sum\limits_{i=0}^h \tau_i=T/2$. 
Defining 
$$
\sigma_0=0,\qquad \sigma_i=\tau_0+...+\tau_{i-1},\qquad i=1,...,h+1.
$$
we have $\sigma_{h+1}=T/2$.
Let $\lambda : [0,T/2]\rightarrow [0,1]$ be the path defined by
\setlength{\extrarowheight}{0.5cm}
\begin{displaymath}
\lambda(t)=\left\{
\begin{array}{lll}
\mu_1\left( 1-\left(\frac{\tau_0-t}{\tau_0}\right)^\frac{2}{3}\right), \quad & \text{if}\quad t\in [0,\sigma_1] & \\
\mu_i+\left(\frac{t-\sigma_i}{\tau_i/2}\right)^\frac{2}{3} \frac{\mu_{i+1}-\mu_i}{2}, \quad & \text{if}\quad 
t\in [\sigma_i,\sigma_i+\frac{\tau_i}{2}],\quad &i=1,...,h-1 \\
\mu_{i+1}- \left(\frac{\sigma_{i+1}-t}{\tau_i/2}\right)^\frac{2}{3} \frac{\mu_{i+1}-\mu_i}{2}, \quad & \text{if}\quad
t\in [\sigma_i+\frac{\tau_i}{2},\sigma_{i+1}], \quad &i=1,...,h-1 \\
\mu_h+(1-\mu_h)\left( \frac{t-\sigma_h}{T/2-\sigma_h}\right)^\frac{2}{3}, \quad & \text{if} \quad t\in 
[\sigma_h,T/2]. & 
\end{array}
\right.
\end{displaymath}
\setlength{\extrarowheight}{0cm}
The definition of $\lambda(t)$ in the interval $[\sigma_0,\sigma_1]$ has some meaning 
only if $\sigma_0<\sigma_1$ (i.e. if 
$\tau_0>0$). 
The path
$$
\gamma_x(t)=(1-\lambda(t))x+\lambda(t)R x_0^\prime, \qquad t\in [0,T/2]
$$
connects $x$ to $R x_0^\prime$ in the time $T/2$. 
If $1\le i\le h-1$, the action of the restriction $\gamma_x\left|_{[\sigma_i,\sigma_i+\tau_i/2]}\right.$ is given
by
\setlength{\extrarowheight}{0.5cm}
\begin{displaymath}
\begin{array}{rl}
A_L(\gamma_x\left|_{[\sigma_i,\sigma_i+\frac{\tau_i}{2}]}\right.)&=\frac{\|R x_0^\prime -x\|^2}{2}
\displaystyle\int_{\sigma_i}^{\sigma_i+\frac{\tau_i}{2}} \dot{\lambda}(t)^2\, dt \\
&+\sum\limits_{1\le j<k\le N} m_j m_k \displaystyle\int_{\sigma_i}^{\sigma_i+\frac{\tau_i}{2}} 
\frac{dt}{|(1-\lambda(t))\vec{r}_{jk}+\lambda(t)R \vec{c}_{jk}|}.
\end{array}
\end{displaymath}
\setlength{\extrarowheight}{0cm}
As $t\in [\sigma_i,\sigma_i+\tau_i/2]$ the path $\lambda(t)$ increases 
from $\mu_i$ to $(\mu_i+\mu_{i+1})/2$, 
hence the coefficient $\lambda_{jk}$ that is closest to $\lambda(t)$ is exactly $\mu_i$.  
Using the triangular inequality we find
\setlength{\extrarowheight}{0.3cm}
\begin{displaymath}
\begin{array}{rl}
|(1-\lambda(t))\vec{r}_{jk}+\lambda(t) R \vec{c}_{jk}| &\ge |r_{jk}-\lambda(t) (r_{jk}+R c_{jk})| \\
&=(r_{jk}+R c_{jk})|\lambda_{jk}-\lambda(t)| \\
&\ge R c_{jk} |\lambda_{jk}-\lambda(t)| \\
&\ge R c_{jk} (\lambda(t)-\mu_i),
\end{array}
\end{displaymath}
\setlength{\extrarowheight}{0cm}
for $t\in [\sigma_i,\sigma_i+\tau_i/2]$ and for every $1\le j<k\le N$.
Therefore, since $\|x\|\le R$ and $\|x_0^\prime\|=1$
\setlength{\extrarowheight}{0.4cm}
\begin{displaymath}
\begin{array}{rl}
A_L(\gamma_x\left|_{[\sigma_i,\sigma_i+\frac{\tau_i}{2}]}\right.)&\le 2R^2 \left(\frac{\mu_{i+1}-\mu_i}{2}\right)^2 
\left(\frac{2}{\tau_i}\right)^\frac{4}{3} \displaystyle\int_{\sigma_i}^{\sigma_i+\frac{\tau_i}{2}} 
\frac{4}{9}(t-\sigma_i)^{-\frac{2}{3}} dt \\ 
&+\sum\limits_{1\le j<k\le N}\frac{m_j m_k}{Rc_{jk}} 
\displaystyle\int_{\sigma_i}^{\sigma_i+\frac{\tau_i}{2}}\frac{dt}{\lambda(t)-\mu_i} \\
&=\frac{4R^2}{3} \frac{(\mu_{i+1}-\mu_i)^2}{\tau_i}+\frac{U(x_0^\prime)}{R}\frac{2}{\mu_{i+1}-\mu_i}
\left(\frac{\tau_i}{2}\right)^\frac{2}{3} \displaystyle\int_{\sigma_i}^{\sigma_i+\frac{\tau_i}{2}} 
\frac{dt}{(t-\sigma_i)^\frac{2}{3}} \\
&=\frac{4R^2}{3} \frac{(\mu_{i+1}-\mu_i)^2}{\tau_i}+\frac{3 U(x_0^\prime)}{R} \frac{\tau_i}{\mu_{i+1}-\mu_i}.
\end{array}
\end{displaymath}
\setlength{\extrarowheight}{0cm}
In a similar way we find
\setlength{\extrarowheight}{0.4cm}
\begin{displaymath}
\begin{array}{rl}
A_L(\gamma_x\left|_{[\sigma_i+\frac{\tau_i}{2},\sigma_{i+1}]}\right.)&\le 
\frac{4R^2}{3} \frac{(\mu_{i+1}-\mu_i)^2}{\tau_i}+\frac{3 U(x_0^\prime)}{R} \frac{\tau_i}{\mu_{i+1}-\mu_i}, \qquad i=1,...,h-1 \\
A_L(\gamma_x\left|_{[\sigma_0,\sigma_1]}\right.)&\le 
\frac{8R^2}{3} \frac{\mu_1^2}{\tau_0}+\frac{3 U(x_0^\prime)}{R} \frac{\tau_0}{\mu_1}, \\
A_L(\gamma_x\left|_{[\sigma_h,\sigma_{h+1}]}\right.)&\le 
\frac{8R^2}{3} \frac{(1-\mu_h)^2}{\tau_h}+\frac{3 U(x_0^\prime)}{R} \frac{\tau_h}{1-\mu_h}.
\end{array}
\end{displaymath}
\setlength{\extrarowheight}{0cm}
That gives 
$$
A_L(\gamma_x\left|_{[\sigma_i,\sigma_{i+1}]}\right.)\le 
\frac{8R^2}{3} \frac{(\mu_{i+1}-\mu_i)^2}{\tau_i}+\frac{6 U(x_0^\prime)}{R} \frac{\tau_i}{\mu_{i+1}-\mu_i}, \qquad i=0,...,h,
$$
and by definition of $\tau_i$ 
\setlength{\extrarowheight}{0.4cm}
\begin{equation} \label{Act-gamma-x}
\begin{array}{rl}
A_L(\gamma_x)&\le \frac{16R^2}{3T} \left( \sum\limits_{i=0}^h (\mu_{i+1}-\mu_i)^\frac{3}{2}\right)
\left(\sum\limits_{i=0}^h (\mu_{i+1}-\mu_i)^\frac{1}{2}\right) \\
&+\frac{3 U(x_0^\prime) T}{R} 
\frac{\sum\limits_{i=0}^h (\mu_{i+1}-\mu_i)^\frac{1}{2}}{\sum\limits_{i=0}^h (\mu_{i+1}-\mu_i)^\frac{3}{2}}.
\end{array}
\end{equation}
\setlength{\extrarowheight}{0cm}
By definition of $\mu_i$ we have 
$$
\mu_{i+1}-\mu_i\ge 0,\qquad \sum\limits_{i=0}^h (\mu_{i+1}-\mu_i)=1.
$$
Let us introduce now the functions
$$
\begin{array}{rl}
f_1:\Rg_+^{h+1}\rightarrow \Rg,\qquad f_1(z)&=\sum\limits_{i=0}^h z_i^\frac{3}{2}, \\
f_2:\Rg_+^{h+1}\rightarrow \Rg,\qquad f_2(z)&=\sum\limits_{i=0}^h z_i^\frac{1}{2}, \\
g:\Rg_+^{h+1}\rightarrow \Rg,\qquad g(z)&=\sum\limits_{i=0}^h z_i, 
\end{array}
$$
and study minima and maxima of $f_1$ and $f_2$ under the condition $g(z)=1$. We show by
induction on $h$ that 
\begin{equation} \label{constraints}
\min\limits_{g(z)=1} f_1(z)=\frac{1}{(h+1)^\frac{1}{2}},\qquad \max\limits_{g(z)=1} f_1(z)=1.
\end{equation}
If $h=0$, condition $g(z)=1$ implies $z_0=1$, thus 
$$
\min\limits_{g(z)=1} f_1(z)=1,\qquad \max\limits_{g(z)=1} f_1(z)=1
$$ 
Assuming now the statement is true up to order $h-1$, let us prove it is true at order $h$. 
By Lagrange multiplier theorem, the unique 
interior critical point of $f_1$ under the condition $g(z)=1$ is given by the equations 
$$
\frac{\partial f_1}{\partial z_i}(z)=\lambda \frac{\partial g}{\partial z_i}(z),\quad i=0,...,h,\qquad \lambda\in\Rg,\qquad
g(z)=1,
$$
this gives 
$$
z_i=\frac{1}{h+1},\quad i=0,...,h,\qquad f_1(z)=\frac{1}{(h+1)^\frac{1}{2}}.
$$
The boundary of the simplex $g(z)=1$ is the set of $z=(z_0,...,z_h)$ such that $\sum\limits_{i=0}^h z_i=1$
and $z_i=0$ for at least one indices $i$. 
By inductive hypothesis, the minimum of $f_1(z)$ on the boundary of $g(z)=1$ is 
$1/h^{1/2}$ and the maximum is $1$. Comparing with the value of $f_1$ on the unique interior critical point of $f_1$ 
we find (\ref{constraints}). In a similar way one prove  
\begin{equation} \label{constraints-2}
\min\limits_{g(z)=1} f_2(z)=1,\qquad \max\limits_{g(z)=1} f_2(z)=(h+1)^\frac{1}{2}.
\end{equation}  
Replacing these estimates in (\ref{Act-gamma-x}) we find
$$
A_L(\gamma_x)\le \frac{16 R^2}{3T}(h+1)^\frac{1}{2}+\frac{3 U(x_0^\prime)T}{R} (h+1),
$$ 
since $h\le \frac{N(N-1)}{2}$, inequality (\ref{PMA}) is proved. 
\vspace{1.5mm} \\ \noindent
We give now the proof of Theorem \ref{thm-exist-solut}. 
\\ \noindent
{\bf Proof of Theorem \ref{thm-exist-solut}}. 
By Propositions \ref{lacapelle-marival} and \ref{cap-de-l-homy}, 
there exist three constants $a>0$, $ \overline{T}>0$ and $\overline{s}>1$ such that
for every $T\ge \overline{T}$, for every $t\ge \overline{s} T$ and for every minimizer
$\overline{\gamma}\in \Sigma(x_i,\gamma_0(t);t)$ we have
\begin{equation} \label{lacapelle-aubareil}
A_L(\overline{\gamma}\left|_{[0,T]}\right.)\le a T^\frac{1}{3}
\end{equation}
Let us prove the equicontinuity of the family 
\begin{equation} \label{cevennes}
\left\{ \overline{\gamma}\left|_{[0,T]}\right.,\quad \overline{\gamma}\in \Sigma(x_i,\gamma_0(t);t),
\quad t\ge \overline{s}T \right\}. 
\end{equation} 
By (\ref{lacapelle-aubareil}) we have 
$$
\int_0^T \|\dot{\overline{\gamma}}(s)\|^2 ds \le 2 a T^\frac{1}{3},
$$
hence, by Cauchy-Schwarz inequality, for every $0\le s<s^\prime\le T$ we have
$$
\begin{array}{rl}
|\overline{\gamma}(s^\prime)-\overline{\gamma}(s)|\le \displaystyle\int_s^{s^\prime} \|\dot{\overline{\gamma}}(u)\| du 
&\le \sqrt{s-s^\prime} 
\left(\displaystyle\int_s^{s^\prime} \|\dot{\overline{\gamma}}(u)\|^2 du \right)^\frac{1}{2} \\
&\le (2a T^\frac{1}{3})^\frac{1}{2} \sqrt{s-s^\prime}.
\end{array}
$$
This gives the equicontinuity of the family (\ref{cevennes}). 
By the way, since $\overline{\gamma}(0)=x_i$, the family 
is also equibounded. By Ascoli theorem we can find a divergent sequence $(t_n)_{n=1}^{+\infty}$  
satisfying $t_n\ge \overline{s} T$ and and a sequence of minimizers
$\gamma_n\in \Sigma(x_i,\gamma_0(t_n);t_n)$ such that the restriction 
$(\gamma_n\left|_{[0,T]}\right.)_{n=1}^{+\infty}$
converges uniformly. 
Applying this argument on an increasing and divergent sequence $(T_k)_{k=1}^{+\infty}$, 
by a diagonal trick we can find an increasing and divergent sequence of times
$(t_n)_{n=1}^{+\infty}$, a sequence of 
minimizers $\gamma_n\in \Sigma(x_i,\gamma_0(t_n);t_n)$ and a path $\gamma : [0,+\infty)\rightarrow \Cac$
such that $(\gamma_{n})_{n=1}^{+\infty}$ converges uniformly to $\gamma$ on every compact interval. 
Moreover, by lower semi-continuity of the action we have 
\begin{equation} \label{la-roque-saint-cristophe}
A_L(\gamma\left|_{[0,T]}\right.)\le \liminf\limits_{n\rightarrow +\infty}
A_L(\gamma_n\left|_{[0,T]}\right.) \le aT^\frac{1}{3}
\end{equation}
for every $T>0$, proving in particular that $A_L(\gamma\left|_{[0,T]}\right.)$ is finite. Therefore, $\gamma(T)$ 
is a non-collision configuration for almost all $T>0$. We prove now that $\gamma$ is a minimizing path. 
Since we want to show that $\gamma\left|_{[0,T]}\right.$ is a minimizer for every $T>0$, it is sufficient to
prove that $\gamma\left|_{[0,T]}\right.$ is a minimizer for $T$ arbitrary great.   
\begin{figure}
\centering
\includegraphics[scale=1.0]{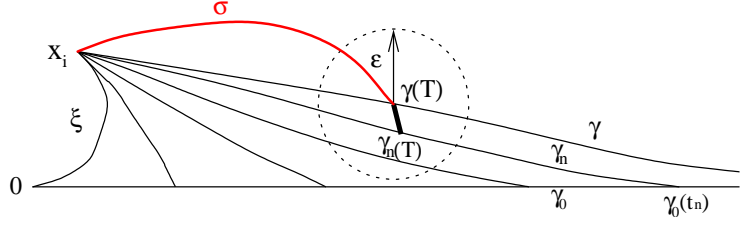}
\caption{$\sigma_{\epsilon,n}$ is obtained by pasting $\sigma$ (reparametrized) with the straight line joining $\gamma(T)$
to $\gamma_n(T)$}
\label{theorem6}
\end{figure}
We can assume, without loss of generality, that $\gamma(T)$ is a non-collision configuration.  
Assuming, for the sake of a contradiction, that $\gamma\left|_{[0,T]}\right.$ is not a minimizer, there would exists an 
absolutely continuous path $\sigma :[0,T]\rightarrow \Cac$ joining $x_i$ to $\gamma(T)$ such that
\begin{equation} \label{la-roque-geac}
A_L(\sigma)<A_L(\gamma\left|_{[0,T]}\right.).
\end{equation}   
Moreover, there exists $M>0$ and $\overline{\epsilon}>0$ such that 
$$
\forall x\in \overline{B}(\gamma(T),\overline{\epsilon}) \Rightarrow U(x)\le M,
$$
where $\overline{B}(\gamma(T),\overline{\epsilon})$ is the closed ball centered in 
$\gamma(T)$ with radius $\overline{\epsilon}$. 
Since the sequence $\gamma_{n}\left|_{[0,T]}\right.$ converges uniformly to 
$\gamma\left|_{[0,T]}\right.$, given $0<\epsilon<\overline{\epsilon}$ 
there exists a positive integer $N_{T,\epsilon}$
such that for every $n\ge N_{T,\epsilon}$ we have $\gamma_{n}(T)\in \overline{B}(\gamma(T),\epsilon)$. 
Let $\sigma_{\epsilon,n}:[0,T]\rightarrow \Cac$ be the path defined by
\setlength{\extrarowheight}{0.4cm}
\begin{displaymath}
\sigma_{\epsilon,n}(t)=\left\{
\begin{array}{ll}
\sigma(\frac{T}{T-\epsilon}t)\qquad &\text{if}\quad t\in [0,T-\epsilon] \\
\frac{T-t}{\epsilon} \gamma(T)+\frac{t-T+\epsilon}{\epsilon} \gamma_{n}(T) \qquad &\text{if}\quad t\in [T-\epsilon,T],
\end{array}
\right.
\end{displaymath}
\setlength{\extrarowheight}{0cm}
where $n\ge N_{T,\epsilon}$.
By construction $\sigma_{\epsilon,n}$ joins $x_i$ to $\gamma_{n}(T)$ in time $T$ (see Figure \ref{theorem6}). 
Moreover, if $t\in [T-\epsilon,T]$, $\sigma_{\epsilon,n}(t)$ is contained in the ball $\overline{B}(\gamma(T),\epsilon)$. 
Computing the action of $\sigma_{\epsilon,n}$ we get 
\setlength{\extrarowheight}{0.4cm}
\begin{displaymath}
\begin{array}{rl}
A_L(\sigma_{\epsilon,n})&\le \frac{1}{2} \left( \frac{T}{T-\epsilon}\right)^2 
\displaystyle\int_0^{T-\epsilon} \left\| \dot{\sigma}\left( \frac{Tt}{T-\epsilon}\right)\right\|^2 dt+
\displaystyle\int_0^{T-\epsilon} U\left(\sigma\left(\frac{Tt}{T-\epsilon}\right)\right)dt \\
&+\left(M+\frac{1}{2}\right)\epsilon  \\
&=\frac{T}{T-\epsilon} \displaystyle\int_0^T \frac{1}{2} \|\dot{\sigma}(t)\|^2 dt+\frac{T-\epsilon}{T} 
\displaystyle\int_0^T U(\sigma(t))dt+
\left(M+\frac{1}{2}\right)\epsilon  \\
&=A_L(\sigma)+{\mathcal O}(\epsilon).
\end{array}
\end{displaymath}
\setlength{\extrarowheight}{0cm}
Inequalities (\ref{la-roque-saint-cristophe}) and (\ref{la-roque-geac}) imply 
$$
A_L(\sigma_{\epsilon, n})<A_L(\gamma_{n}\left|_{[0,T]}\right.)
$$
if $\epsilon$ is sufficiently small and $n$ sufficiently great. This contradicts the minimizing property of $\gamma_{n}$ 
and proves that $\gamma$ is a minimizer. By Marchal theorem, 
$\gamma$ is collision-free (and in particular it is a real solution of the
N-body problem) for $t>0$. This complete the proof of Theorem \ref{thm-exist-solut}.  
%
%
%
%
%
%
%
%
%
%
\vspace{2mm}
\section{\bf Parabolicity of the solution} \label{parabol-solut}
\vspace{3mm} \noindent
To complete the proof of the main theorem we still have to show that the limit solution $\gamma(t)$ 
is parabolic and asymptotic to $x_0$. By Lemma \ref{lem-parab-asympt-x0} we just need to
verify the asymptotic estimates (\ref{gamma-I_0-t43}).
We introduce now the following 
\begin{notat} 
Given the functions $f(r,x_1,...,x_n)$ and $g(r,x_1,...,x_n)\neq 0$, we write 
$f(r,x_1,...,x_n)=o_r(g(r,x_1,...,x_n))$ as $r\rightarrow r_0$ if the quotient 
$\frac{f(r,x_1,...,x_n)}{g(r,x_1,...,x_n)}$ is infinitesimal as $r\rightarrow r_0$, uniformly on $(x_1,...,x_n)$. 
In a similar way, we write $f(r,x_1,...,x_n)={\mathcal O}_r(g(r,x_1,...,x_n))$ 
if the quotient $\frac{f(r,x_1,...,x_n)}{g(r,x_1,...,x_n)}$ is locally bounded for $r$ close to $r_0$, 
uniformly on the variables $(x_1,...,x_n)$.
\end{notat}
\noindent
Let us give now a refinement of Lemma \ref{lemme-fondamentale}.
\begin{lem} \label{nyons}
Let $\tau>0$ and $\xi\in\Sigma(0,x_i;\tau)$.  
There exist two constants $\overline{T}>\tau$ and $\overline{s}>1$ such that for every $T\ge \overline{T}$, 
for every $t\ge \overline{s} T$ and for every minimizer $\overline{\gamma}\in \Sigma(x_i,\gamma_0(t);t)$
we have 
\setlength{\extrarowheight}{0.3cm}
\begin{displaymath}
\begin{array}{rl}
{\mathcal M}(\overline{\gamma}(T),\gamma_0(t);0,T,t)&\le 2A_L(\xi)+{\mathcal O}_T(T^{-\frac{2}{3}}) \\
{\mathcal N}(\|\overline{\gamma}(T)\|,\alpha t^\frac{2}{3};0,T,t)&\le 2A_L(\xi)+{\mathcal O}_T(T^{-\frac{2}{3}}). 
\end{array}
\end{displaymath}
\setlength{\extrarowheight}{0cm}
as $T\rightarrow +\infty$.
\end{lem}
\Pf\ 
The second inequality is a direct consequence of the first one and of Lemma \ref{nonna-isabella}. 
Let us prove the first inequality.
We consider $\tau$ as a fixed constant, while $T$ and $t$ are variables. Let $\overline{T}>0$, $\overline{s}>1$ and $K>0$ be like in Proposition \ref{lacapelle-marival}. 
Without loss of generality we can assume $\overline{T}>\tau$.  
Let $T\ge \overline{T}$ and $t\ge \overline{s}T$. Let $\eta_{T+\tau}\in\Sigma(0,\overline{\gamma}(T);T+\tau)$
The path 
$$
\overline{\eta}_T: [0,T]\rightarrow \Cac,\qquad \overline{\eta}_T(s)=\eta_{T+\tau}\left(\frac{T+\tau}{T}s\right).
$$
is a reparametrization of $\eta_{T+\tau}$ and it joins $0$ to $\overline{\gamma}(T)$ in time $T$, thus 
$$
{\mathcal A}(0,\overline{\gamma}(T);T)\le A_L(\overline{\eta}_T).
$$
A computation of the action of $\overline{\eta}_T$ gives 
\begin{equation} \label{taybosc}
A_L(\overline{\eta}_T)=\left(1+{\mathcal O}(1/T)\right) A_L(\eta_{T+\tau}), \qquad T\rightarrow +\infty.
\end{equation}
Since $\eta_{T+\tau}$ is a minimizer joining $0$ to $\overline{\gamma}(T)$ 
in time $T+\tau$, 
by Propositions \ref{lacapelle-marival} and \ref{cap-de-l-homy} we obtain 
\begin{equation} \label{cazaubon}
A_L(\eta_{T+\tau})\le C_1 \frac{K^2 T^\frac{4}{3}}{T+\tau}+C_2 \frac{T+\tau}{K T^\frac{2}{3}}={\mathcal O}(T^\frac{1}{3}).
\end{equation}
Combining inequalities (\ref{taybosc}) and  (\ref{cazaubon}), by definition of $\eta_{T+\tau}$ and $\overline{\eta}_T$
we get
\begin{equation} \label{barbontan}
{\mathcal A}(0,\overline{\gamma}(T);T)-{\mathcal A}(0,\overline{\gamma}(T);T+\tau)\le {\mathcal O}(T^{-\frac{2}{3}}). 
\quad T\rightarrow +\infty
\end{equation}
In a similar way, let us consider a minimizer $\eta_T\in\Sigma(0,\overline{\gamma}(T);T)$ 
The path 
$$
\overline{\eta}_{T+\tau}:[0,T+\tau]\rightarrow \Cac,\qquad \overline{\eta}_{T+\tau}(s)=\eta_T\left(\left(\frac{T}{T+\tau}\right)s\right).
$$ 
is a reparametrization of $\eta_T$, and it joins $0$ to $\overline{\gamma}(T)$ in time $T+\tau$, hence
$$
{\mathcal A}(0,\overline{\gamma}(T);T+\tau)\le A_L(\overline{\eta}_{T+\tau}).
$$
Arguing as before we get the estimates
\begin{equation} \label{pech-merle}
{\mathcal A}(0,\overline{\gamma}(T);T+\tau)-{\mathcal A}(0,\overline{\gamma}(T);T)\le {\mathcal O}(T^{-\frac{2}{3}}), 
\quad T\rightarrow +\infty.
\end{equation}
Combining inequalities (\ref{barbontan}) with (\ref{pech-merle}) we obtain
\begin{equation} \label{cro-magnon}
{\mathcal A}(0,\overline{\gamma}(T);T+\tau)={\mathcal A}(0,\overline{\gamma}(T);T)+ {\mathcal O}_T(T^{-\frac{2}{3}}), 
\quad T\rightarrow +\infty,
\end{equation}
uniformly on $t\ge \overline{s} T$ and $\overline{\gamma}\in\Sigma(0,\gamma_0(t);t)$. 
With the same argument we find the following estimates
\begin{equation} \label{dentelles}
{\mathcal A}(0,\gamma_0(t);t-\tau)-{\mathcal A}(0,\gamma_0(t);t)={\mathcal O}_t(t^{-\frac{2}{3}}), 
\quad t\rightarrow +\infty.
\end{equation}
Replacing (\ref{cro-magnon}) and (\ref{dentelles}) into the first inequality of 
Lemma \ref{lemme-fondamentale}, since we assume $t\ge \overline{s}T$ and $\overline{s}>1$, we obtain the first inequality
of this Lemma. This ends the proof. 
\vspace{0.5mm} \\ \noindent
To simplify the notations we introduce now the functions 
\setlength{\extrarowheight}{0.2cm}
\begin{displaymath}
\begin{array}{rl}
&{\mathcal F}:\Cac\times (1,+\infty)\rightarrow \Rg_+ \\
&{\mathcal F}(x,s)={\mathcal M}(x,\gamma_0(s);0,1,s)={\mathcal A}(0,x;1)+{\mathcal A}(x,\gamma_0(s);s-1)
-{\mathcal A}(0,\gamma_0(s);s),
\end{array}
\end{displaymath}
and 
\begin{displaymath}
\begin{array}{rl}
&{\mathcal G}:\Rg_+\times (1,+\infty)\rightarrow \Rg_+ \\
&{\mathcal G}(r,s)={\mathcal N}(r,\alpha s^\frac{2}{3};0,1,s)=S(0,r;1)+S(r,\alpha s^\frac{2}{3};s-1)
-S(0,\alpha s^\frac{2}{3};s).
\end{array}
\end{displaymath}
\setlength{\extrarowheight}{0cm} 
\begin{lem} \label{sommieres}
Given $s>1$ and $x\in\Cac$ we have 
$$
{\mathcal F}(x,s)\ge 0
$$
with equality if and only if $x=\alpha x_0$.
\end{lem}
\Pf\ 
By Lemma \ref{nonna-isabella} we have 
$$
{\mathcal F}(x,s)\ge {\mathcal G}(\|x\|,s),
$$
with equality if and only if $x=\|x\| x_0$. Since $u\mapsto \alpha u^\frac{2}{3}$ is the unique
solution of the one-dimensional Kepler problem joining $0$ to $\alpha s^\frac{2}{3}$ in time $s$ (see Lemma 
\ref{uniq-solution} in the Appendix), it is also a
minimizer, therefore 
$$
{\mathcal G}(r,s)\ge 0,
$$ 
with equality if and only if $r=\alpha$. This proves the Lemma.
\vspace{0.5mm} \\ \noindent 
By homothety invariance, the conclusion of Lemmas \ref{nyons} and \ref{sommieres} can be written in the more 
compact form 
\setlength{\extrarowheight}{0.3cm}
\begin{equation} \label{nesque} 
\begin{array}{rl}
0\le T^{1/3}{\mathcal F}\left(\frac{\overline{\gamma}(T)}{T^{2/3}},\frac{t}{T}\right)&\le 2A_L(\xi)
+{\mathcal O}_T(T^{-{2/3}}) \\
0\le T^{1/3}{\mathcal G}\left(\frac{\|\overline{\gamma}(T)\|}{T^{2/3}},\frac{t}{T}\right)&\le 2A_L(\xi)
+{\mathcal O}_T(T^{-{2/3}}).
\end{array}
\end{equation}
\setlength{\extrarowheight}{0cm}
as $T\rightarrow +\infty$, uniformly on $t\ge \overline{s}T$ and $\overline{\gamma}\in\Sigma(0,\gamma_0(t);t)$.  \\
The following Theorem is a main tool in the proof of the Main Theorem. It shows that if 
${\mathcal F}(x,s)$ is sufficiently small and $s$ is sufficiently great, the configuration $x$ is 
close to $\alpha x_0$. 
\begin{thm} \label{ardeche}
There exist a function $\delta : (0,\overline{\epsilon}]\rightarrow \Rg_+$
satisfying $\delta(\epsilon)=o(1)$ 
as $\epsilon\rightarrow 0^+$, such that for every $\epsilon\in (0,\overline{\epsilon}]$, 
there exists $\overline{s}_\epsilon>1$, such that for  
every $s\ge \overline{s}_\epsilon$, the set of configurations $x\in\Cac$ 
satisfying the inequality 
\begin{equation} \label{annie}
{\mathcal F}(x,s)\le \epsilon
\end{equation}
is contained in the ball ${\overline B}(\alpha x_0,\delta(\epsilon))$. 
\end{thm}
\noindent  
Before giving the proof of Theorem \ref{ardeche}, 
we show that this theorem achieve the proof of the Main Theorem.
\\ \noindent
{\bf Proof of the Main Theorem}.  Let $\gamma :[0,+\infty)\rightarrow \Cac$ be the limit solution 
constructed in Theorem \ref{thm-exist-solut} and let $\gamma_n\in \Sigma(x_i,\gamma_0(t_n);t_n)$ be the sequence of 
minimizers uniformly convergent to $\gamma$ on every compact interval. Let $\overline{\epsilon}$ like in Theorem \ref{ardeche}
and $0<\epsilon<\overline{\epsilon}$. An immediate consequence of inequalities (\ref{nesque}) is the existence of 
$\overline{T}_\epsilon\ge \overline{T}$ such that if $T\ge \overline{T}_\epsilon$ and $t_n\ge \overline{s} T$ we have 
$$
{\mathcal F}\left( \frac{\gamma_n(T)}{T^\frac{2}{3}},\frac{t_n}{T}\right)\le \epsilon
$$ 
and by Theorem \ref{ardeche}
$$
\left\|\frac{\gamma_n(T)}{T^{2/3}}-\alpha x_0\right\|\le \delta(\epsilon),
$$
for $t_n$ sufficiently great. 
The sequence $\gamma_n\left|_{[0,T]}\right.$ converges uniformly to $\gamma\left|_{[0,T]}\right.$ 
as $n\rightarrow +\infty$, hence
$$
\left\|\frac{\gamma(T)}{T^{2/3}}-\alpha x_0\right\|\le \delta(\epsilon),
$$
for every $T\ge \overline{T}_\epsilon$. Since $\delta(\epsilon)\rightarrow 0$ as 
$\epsilon\rightarrow 0$, we have proved that
$$
\frac{\gamma(T)}{T^{2/3}}\rightarrow \alpha x_0,\quad \text{as}\quad T\rightarrow +\infty,
$$   
that is to say, $\gamma$ is parabolic and asymptotic to $x_0$. This achieves the proof of the Main Theorem. 
\\ \noindent
The next section is devoted to prove Theorem \ref{ardeche}.
%
%
%
%
\vspace{2mm}
\section{\bf Proof of Theorem \ref{ardeche}}   \label{Lambert}
\vspace{3mm} \noindent
In order to achieve the proof of Theorem \ref{ardeche} we compare the $N$-body problem with a Kepler problem 
on the configuration space with a lagrangian given by   
$$
L_0(x,\dot{x})=\frac{\|\dot{x}\|^2}{2}+\frac{U_0}{\|x\|},\qquad (x,\dot{x})\in\Cac\times \Cac. 
$$
Let $A_{L_0}(\varpi)$ denote the action (for the lagrangian $L_0$) of an absolutely continuous path $\varpi$ and 
${\mathcal A}_0(x_1,x_2;s)$ the infimum of $A_{L_0}(\varpi)$ over all absolutely continuous paths $\varpi$  
joining $x_1$ to $x_2$ in time $s$. We have the inequality 
$$
{\mathcal A}(x_1,x_2;s) \ge {\mathcal A}_0(x_1,x_2;s)\ge S(\|x_1\|,\|x_2\|;s),
$$
with ${\mathcal A}(x_1,x_2;s)={\mathcal A}_0(x_1,x_2;s)$  
if and only if there exists a minimizing path (for the lagrangian $L$) $\varpi :[0,s]\rightarrow \Cac$ 
joining $x_1$ with $x_2$ such that $\tilde{U}(\varpi(u))=U_0$ for every $u\in [0,s]$, and 
${\mathcal A}_0(x_1,x_2;s)=S(\|x_1\|,\|x_2\|;s)$ if and only if $x_1$ and $x_2$ are on a same half-line starting 
from the origin. 
The function     
\setlength{\extrarowheight}{0.2cm}
\begin{displaymath}
\begin{array}{rl}
&{\mathcal F}_0:\Cac\times (1,+\infty)\rightarrow \Rg_+, \\
&{\mathcal F}_0(x,s)={\mathcal A}_0(0,x;1)+{\mathcal A}_0(x,\gamma_0(s);s-1)-{\mathcal A}_0(0,\gamma_0(s);s),
\end{array}
\end{displaymath}
\setlength{\extrarowheight}{0cm}
verifies the inequality
\begin{equation} \label{pas-de-la-graille}
{\mathcal F}(x,s)\ge {\mathcal F}_0(x,s)\ge {\mathcal G}(\|x\|,s)\ge 0.
\end{equation}
Roughly speaking, to achieve the proof of Theorem \ref{ardeche}, we replace ${\mathcal F}(x,s)$ 
with ${\mathcal F}_0(x,s)$ 
and we show that if $\epsilon$ is small and $s$ great, the inequality ${\mathcal F}_0(x,s)\le \epsilon$ can be satisfied 
only if $x$ is in a small ball centered in $\alpha x_0$. \\ \noindent
This goal will be achieved in two steps. In Proposition \ref{Gcors(r)} we prove that if $s$ is sufficiently great, 
the set of $r\in\Rg_+$ verifying ${\mathcal G}(r,s)\le \epsilon$ is contained in a small interval centered in 
$\alpha$. Hence, by inequality (\ref{pas-de-la-graille}), the set of configuration $x$ verifying 
${\mathcal F}(x,s)\le \epsilon$ is contained in a thin hollow sphere with inner and outer radious close to 
$\alpha$.
In Proposition \ref{applic-lambert} we show that the set of configurations $x$ 
verifying ${\mathcal F}_0(x,s)\le \epsilon$ is a small neighborhood of $\alpha x_0$.
\begin{prop} \label{Gcors(r)} 
There exist a function $\delta_1 : (0,\overline{\epsilon}_1]\rightarrow \Rg_+$ satisfying $\delta_1(\epsilon)=o(1)$ as
$\epsilon\rightarrow 0^+$, such that for every $\epsilon\in (0,\overline{\epsilon}_1]$ there exists 
$\overline{s}^1_\epsilon>1$, such that for every $s\ge \overline{s}^1_\epsilon$, 
the set of $r\in\Rg_+$ satisfying the inequality 
$$
{\mathcal G}(r,s)\le \epsilon
$$
is contained in the interval $[\alpha-\delta_1(\epsilon),\alpha+\delta_1(\epsilon)]$ . 
\end{prop}  
\Pf\ By Proposition \ref{elenuccia} of the Appendix there exists $\overline{r}>0$ 
and $\overline{s}>0$ such that for every $r\ge \overline{r}$ and for every $s\ge \overline{s}$ 
we have ${\mathcal G}(r,s)> 1$. Without
loss of generality we will assume $\alpha<\overline{r}<\overline{s}^{1/3}$.  
By Proposition \ref{S(c-point)} of the Appendix we have 
$$
{\mathcal G}(r,s)=S(0,r;1)-\beta_0 r^{\frac{1}{2}}+g(r,s), 
$$
where $g(r,s)=o_s(1)$ as $s\rightarrow +\infty$, uniformly on $0\le r \le s^\frac{1}{3}$, and where $\beta_0=(8 U_0)^\frac{1}{2}$.
Let us introduce now the function
$$
G(r)=S(0,r;1)-\beta_0 r^\frac{1}{2}.
$$
By Lemma \ref{uniq-solution} the solution joining $0$ with $r$
in time $1$ is monotonic if and only if for $r\ge \beta$, where 
$\beta=2\left(\frac{U_0}{\pi^2}\right)^{1/3}$. We remark that $\beta<\alpha$. 
The energy $h(0,r;1)$ of this solution 
is negative if and only if $0\le r<\alpha$, moreover $h(0,\beta;1)=-U_0/\beta$.  
Let us term $h=h(0,r;1)$. The action $S(0,r;1)$ is given by
\setlength{\extrarowheight}{0.6cm}
\begin{displaymath}
S(0,r;1)=\left\{
\begin{array}{ll}
\displaystyle\int_0^{-\frac{U_0}{h}} \sqrt{2\left(h+\frac{U_0}{u} \right)}du
+\displaystyle\int_r^{-\frac{U_0}{h}} \sqrt{2\left(h+\frac{U_0}{u} \right)}du -h\quad &\text{if}\quad r <\beta \\
\displaystyle\int_0^r \sqrt{2\left(h+\frac{U_0}{u} \right)}du -h \quad &\text{if}\quad r\ge \beta,
\end{array}
\right.
\end{displaymath}
\setlength{\extrarowheight}{0cm}
hence by Lemma \ref{lem-h-croiss}, the functions $r\mapsto S(0,r;1)$ and $r\mapsto G(r)$ are of class $\Cr^1$ 
on $(0,+\infty)$,  moreover we have 
\setlength{\extrarowheight}{0.4cm}
\begin{displaymath}
G^\prime(r)=\left\{
\begin{array}{ll}
-\sqrt{2\left(h(0,r;1)+\frac{U_0}{r}\right)}-\sqrt{\frac{2U_0}{r}}&\text{if}\quad 0<r<\beta \\
\sqrt{2\left(h(0,r;1)+\frac{U_0}{r}\right)}-\sqrt{\frac{2U_0}{r}} \quad &\text{if}\quad r\ge \beta,
\end{array}
\right.
\end{displaymath}
\setlength{\extrarowheight}{0cm}
proving that $G(r)$ is in fact of class $\Cr^2$ on $(0,\beta)\cup (\beta,+\infty)$.
Since the function $r\mapsto h(0,r;1)$ is increasing and $h(0,\alpha;1)=0$, the function
$G(r)$ is decreasing for $r\in (0,\alpha)$ and it is increasing for $r\in (\alpha,+\infty )$. 
The absolute minimum of $G(r)$ is achieved at $r=\alpha$, and we have 
$$
G(\alpha)=\int_0^{\alpha}\sqrt{\frac{2U_0}{u}} du-\sqrt{8 U_0 \alpha}=0.
$$  
By Lemma \ref{lem-h-croiss}, a direct computation of the second derivative of $G$ at $\alpha$ gives
$$
G^{\prime\prime}(\alpha)=\frac{5 U_0^\frac{1}{2}}{2^\frac{1}{2} \alpha^\frac{3}{2}},
$$
hence, since $G(\alpha)=G^\prime(\alpha)=0$, there exists 
$\overline{\delta}>0$ and $C_1>0$ such that 
$$
\forall r\in [\alpha-\overline{\delta},\alpha+\overline{\delta}],\qquad G(r)\ge C_1(r-\alpha)^2. 
$$
Without loss of generality we shall assume $\alpha-\overline{\delta}>\beta$ and 
$\alpha+\overline{\delta}<\overline{r}$. 
Let $\overline{\epsilon}_1=\min\{\frac{C_1 \overline{\delta}^2}{2}, 1\}$ and let us define the function
$$
\delta_1 :(0,\overline{\epsilon}_1]\rightarrow \Rg_+,\qquad \delta_1(\epsilon)=\sqrt{\frac{2\epsilon}{C_1}}.
$$
Since $G(r)$ is decreasing for $r\le \alpha$ and increasing for $r\ge \alpha$,  
for every $\epsilon\in (0,\overline{\epsilon}_1]$ we have 
\begin{equation} \label{G-r} 
\forall r\in (0,\alpha-\delta_1(\epsilon))\cup (\alpha+\delta_1(\epsilon),+\infty),\qquad 
G(r)> C_1 \delta_1(\epsilon)^2=2\epsilon.
\end{equation}
We come back now to the function ${\mathcal G}(r,s)=G(r)+g(r,s)$. Since $g(r,s)$ is infinitesimal
for $s\rightarrow +\infty$ and $0\le r\le s^\frac{1}{3}$, for every $\epsilon\in (0,\overline{\epsilon}_1]$ 
there exists $\overline{s}^1_\epsilon>\overline{s}$ such that for every $s\ge \overline{s}^1_\epsilon$ 
and for every $r$ verifying $0\le r\le s^\frac{1}{3}$ we have $|g(r,s)|\le \epsilon$. 
If $s\ge \overline{s}^1_\epsilon$ and $r\ge \overline{r}$ we have 
$$
{\mathcal G}(r,s)>1\ge \overline{\epsilon}_1\ge \epsilon.
$$    
If $s\ge \overline{s}^1_\epsilon$ and $r\in ]0, \overline{r}[$, by (\ref{G-r}), for every 
$r\in ]0,\alpha-\delta_1(\epsilon)[\cup ]\alpha+\delta_1(\epsilon),\overline{r}[$ we have 
$$
{\mathcal G}(r,s)=G(r)+g(r,s)>2\epsilon-\epsilon=\epsilon.
$$ 
This ends the proof of the Proposition.
\vspace{0.5mm}\\ \noindent
We introduce the following notation : given two configurations $x_1$ and $x_2$, the angle between $x_1$ and $x_2$
is denoted by the symbol $\angle(x_1,x_2)$. We always have $0\le \angle(x_1,x_2) \le \pi$.  
%
%
%
%
%
%
\begin{prop} \label{applic-lambert}
If $\overline{\epsilon}_1$ and $\delta_1 :(0,\overline{\epsilon}_1]\rightarrow \Rg_+$ are like in 
Proposition \ref{Gcors(r)}, there exist $\overline{\epsilon}\in (0,\overline{\epsilon}_1]$ and $C_2>0$ 
such that given the function 
\begin{equation} \label{delta_2}
\delta_2 :(0,\overline{\epsilon}]\rightarrow \Rg_+, \qquad \delta_2(\epsilon)=(C_2\epsilon)^\frac{1}{2},
\end{equation}
for every $\epsilon\in (0,\overline{\epsilon}]$, there exists $\overline{s}^2_\epsilon>1$ such that for  
every $s\ge \overline{s}^2_\epsilon$ and for every configuration $x\in\Cac$ satisfying
\begin{equation} \label{angle-config}
\left|\|x\|-\alpha\right|\le \delta_1(\epsilon),\qquad \angle(x,x_0)> \delta_2(\epsilon)
\end{equation}
we have
$$
{\mathcal F}_0(x,s) >\epsilon.
$$
\end{prop}
\Pf\ 
The basic tool of this proof is Lambert's Theorem. Our reference is [Al]. \\ \noindent
Let $C_2>0$ and $\overline{\epsilon}\in (0,\overline{\epsilon}_1]$. 
Let $\delta_2 : (0,\overline{\epsilon}]\rightarrow \Rg_+$ be the function defined in (\ref{delta_2}). 
In the following we will ask more precise conditions on $C_2$ and $\overline{\epsilon}$.
Let $0<\epsilon\le \overline{\epsilon}$, let $x$ be a configuration verifying (\ref{angle-config}) and
$s>1$. 
The minimizer (for $L_0$) $\sigma : [0,s-1]\rightarrow \Cac$ joining $x$ to $\gamma_0(s)$ 
in time $s-1$ is a collision-free Keplerian arc, hence it is contained in the plane generated by $0$, 
$x$ and $\gamma_0(s)$. 
Introducing a system of polar coordinates in this plane, we can identify $x$ with $r e^{\imath \theta}$ and    
$\gamma_0(s)$ with $\alpha s^\frac{2}{3}\in \Rg\subset \Cg$ 
where 
$$
|r-\alpha|\le \delta_1(\epsilon),\qquad \delta_2(\epsilon)< |\theta| \le \pi.
$$
Moreover, the path $\sigma$ can be written in polar coordinates by 
$$
\sigma(u)= \rho(u) e^{\imath \phi(u)},\quad u\in [0,s-1],
$$
where 
$$
\begin{array}{rl}
\rho(0)=r \qquad & \phi(0)=\theta \\
\rho(s-1)=\alpha s^\frac{2}{3} \qquad & \phi(s-1)\in 2\pi \Zg.
\end{array}
$$
Since $\sigma$ is collision-free, $\rho(u)>0$ for all $u\in [0,s-1]$. 
By definition of ${\mathcal F}_0$ and using the properties of ${\mathcal A}_0$ we have 
$$
\begin{array}{rl}
{\mathcal F}_0(r e^{\imath \theta},s)&={\mathcal A}_0(0,r e^{\imath \theta};1)+{\mathcal A}_0(r e^{\imath \theta},\gamma_0(s);s-1)
-{\mathcal A}_0(0,\gamma_0(s);s) \\
&=S(0,r;1)+{\mathcal A}_0(r e^{\imath \theta},\gamma_0(s);s-1)-S(0,\alpha s^\frac{2}{3};s).
\end{array}
$$
We prove now that $\sigma$ is a {\sl direct path}, that is to say, the total variation of the polar angle
$\phi$ is less or equal to $\pi$. Assume, for the sake of contradiction, that $|\phi(s-1)-\phi(0)|>\pi$. Eventually 
changing the orientation of the plane, we can assume without loss of generality 
$\phi(s-1)-\phi(0)>\pi$, hence there exists a unique integer $k\ge 1$ and a unique real number
$\alpha\in (-\pi,\pi]$ such that 
$$
\phi(s-1)-\phi(0)=2k\pi+\alpha.
$$    
The path $\overline{\rho} e^{\imath \overline{\phi}}$ defined by 
$$
\overline{\rho}(u)=\rho(u),\quad \overline{\phi}(u)=\phi(0)+\frac{\alpha}{2k\pi+\alpha} (\phi(u)-\phi(0)),
$$
has the same ends as the original one, moreover 
$$
A_{L_0}(\overline{\rho} e^{\imath \overline{\phi}})-A_{L_0}(\rho e^{\imath \phi})=
\frac{1}{2}\left[\left(\frac{\alpha}{2k\pi+\alpha}\right)^2-1\right]\int_0^{s-1} (\rho^2\dot{\phi}^2)(u)du<0,
$$
and we get a contradiction. Lambert's Theorem state that if $x_1$ and $x_2$ are
two configurations and $\tau>0$, the action ${\mathcal A}_0(x_1,x_2;\tau)$ of the direct
Keplerian arc joining $x_1$ to $x_2$ in time $\tau$ is a function of three parameters only : the time
$\tau$, the distance $\|x_1-x_2\|$ between the two ends and the sum of the distances 
between the ends and the origin (i.e. $\|x_1\|+\|x_2\|$). 
Comparing now $\sigma$ with a direct collinear arc, by Lambert's Theorem we find  
$$
{\mathcal A}_0(r e^{\imath \theta},\gamma_0(s);s-1)=S(d_1(r,\theta,s),d_2(r,\theta,s);s-1),
$$
where 
\setlength{\extrarowheight}{0.2cm}
\begin{displaymath}
\begin{array}{rl}
d_1(r,\theta,s)&=\frac{r+\alpha s^\frac{2}{3}-|r e^{\imath \theta}-\alpha s^{\frac{2}{3}}|}{2}, \\ 
d_2(r,\theta,s)&=\frac{r+\alpha s^\frac{2}{3}+|r e^{\imath \theta}-\alpha s^{\frac{2}{3}}|}{2}.
\end{array}
\end{displaymath}
\setlength{\extrarowheight}{0cm}
Moreover
$$
|r e^{\imath \theta}-\alpha s^\frac{2}{3}|=\alpha s^\frac{2}{3}-r\cos \theta+l(r,\theta,s),
$$
where 
$$
l(r,\theta,s)={\mathcal O}_s(s^{-{2/3}}),\quad s\rightarrow +\infty
$$ 
uniformly on $\delta_2(\epsilon)<|\theta|\le \pi$ and 
$|r-\alpha|\le \delta_1(\epsilon)$. 
Therefore we get 
\setlength{\extrarowheight}{0.3cm} 
\begin{displaymath}
\begin{array}{rl}
d_1(r,\theta,s)&=r\left(\frac{1+\cos\theta}{2}\right)-\frac{l(r,\theta,s)}{2} \\ 
d_2(r,\theta,s)&=\alpha s^\frac{2}{3}+r\left(\frac{1-\cos\theta}{2}\right)+\frac{l(r,\theta,s)}{2}.
\end{array}
\end{displaymath}
\setlength{\extrarowheight}{0cm}
Since $S(0,\alpha s^\frac{2}{3};s)=\alpha_0 s^\frac{1}{3}$, applying Proposition \ref{S(c-point)} of the Appendix to
$S(d_1(r,\theta,s),d_2(r,\theta,s);s-1)$ we find
$$
{\mathcal F}_0(re^{\imath \theta},s)=G(r)+\beta_0 r^\frac{1}{2}\left[ 1-\left(\frac{1+\cos\theta}{2}
-\frac{l(r,\theta,s)}{2r}\right)^\frac{1}{2}\right]+g(r,\theta,s),
$$
where $g(r,\theta,s)$ is infinitesimal as $s\rightarrow +\infty$, uniformly on $r$ and $\theta$. In Proposition
\ref{Gcors(r)} we showed that $G(r)\ge 0$ for all $r>0$. Let $\overline{s}^2_\epsilon>0$ such that for every
$s\ge\overline{s}^2_\epsilon$, for every $\theta$ satisfying $|\theta|\in (\delta_2(\epsilon),\pi]$ and for every 
$r\in [\alpha-\delta_1(\epsilon),\alpha+\delta_1(\epsilon)]$ we have 
$$
|g(r,\theta,s)|\le \epsilon,\qquad \left|\frac{l(r,\theta,s)}{2r}\right|\le \epsilon. 
$$  
Since the function $x\mapsto \cos x$ is decreasing in $[0,\pi]$, chosing $C_2>4$ and using the classical expansions 
of $\cos x$ and $(1+x)^\frac{1}{2}$ we find
\setlength{\extrarowheight}{0.3cm} 
\begin{displaymath}
\begin{array}{rl}
{\mathcal F}_0(r e^{\imath \theta},s) &\ge \beta_0 (\alpha-\delta_1(\epsilon))^\frac{1}{2}
\left[1-\left(\frac{1+\cos \delta_2(\epsilon)}{2}+\epsilon\right)^\frac{1}{2}\right]-\epsilon \\
&=\epsilon \left[ \beta_0 \alpha^\frac{1}{2} \left( 1-\frac{\delta_1(\epsilon)}{\alpha}\right)^\frac{1}{2} 
\left( \frac{C_2-4}{8}+\mu(\epsilon)\right)-1 \right],
\end{array}
\end{displaymath}
\setlength{\extrarowheight}{0cm} 
where $\mu(\epsilon)=o(1)$ as $\epsilon\rightarrow 0$. 
Chosing $0<\overline{\epsilon}<\overline{\epsilon}_1$ such that 
$$
\forall \epsilon \in (0,\overline{\epsilon}], \qquad |\mu(\epsilon)|<\frac{1}{8},\quad \text{and}\quad 
|\delta_1(\epsilon)|<\frac{\alpha}{2},
$$ 
and chosing $C_2$ in such a way  
$$
C_2>5+\frac{16\sqrt{2}}{\beta_0 \alpha^\frac{1}{2}}
$$     
we find 
$$
{\mathcal F}_0(r e^{\imath \theta},s)> \epsilon,
$$
for every $r\in [\alpha-\delta_1(\epsilon),\alpha+\delta_1(\epsilon)]$, for every $\theta$ such that 
$|\theta| \in (\delta_2(\epsilon),\pi]$ 
and for every $s\ge \overline{s}^2_\epsilon$. This proves the Proposition. 
\\ \noindent
The proof of Theorem \ref{ardeche} is essentially the juxtaposition of the two previous Propositions. 
\\ \noindent
{\bf Proof of Theorem \ref{ardeche}}. 
\\ \noindent
We use the same notations of the previous two Propositions. Given $\epsilon\in (0,\overline{\epsilon}]$, let  
$\overline{s}_\epsilon=\max\{\overline{s}^1_\epsilon, \overline{s}^2_\epsilon \}$. 
By Proposition \ref{Gcors(r)} and \ref{applic-lambert} and by inequality (\ref{pas-de-la-graille}), if 
$s\ge \overline{s}_\epsilon$ and $x$ is a configuration verifying ${\mathcal F}(x,s)\le \epsilon$ 
we have 
\begin{equation} \label{condit-presque-conclus}
\left|\|x\|-\alpha\right|\le \delta_1(\epsilon)\quad \text{and}\quad \angle(x,x_0)\le \delta_2(\epsilon).
\end{equation}
Let $\delta$ be the function
$$
\delta :(0,\overline{\epsilon}]\rightarrow \Rg_+,\qquad \delta(\epsilon)=\left[2\alpha\left(\alpha+\delta_1(\epsilon)\right)
(1-\cos\delta_2(\epsilon))+\delta_1(\epsilon)^2\right]^\frac{1}{2},
$$ 
an easy computation show that $\delta(\epsilon)\rightarrow 0$ as $\epsilon\rightarrow 0$ and the set of configurations verifying
(\ref{condit-presque-conclus}) is contained in the ball $\overline{B}\left(\alpha x_0,\delta(\epsilon)\right)$.
The Theorem is proved. 
%
%
%
%
\vspace{2mm}
\section*{\bf Appendix : Some estimates for the one-dimensional Kepler Problem} \label{estim-kepler}
\vspace{3mm} \noindent
The Kepler problem on the half-line $\Rg_+$ is defined by the equation 
\begin{equation} \label{Kepler-line}
\ddot{r}=-\frac{U_0}{r^2},
\end{equation}
where $U_0>0$ is the gravitational constant. The Lagrangian function of the problem and the energy 
are written 
$$
l=\frac{\dot{r}^2}{2}+\frac{U_0}{r},\qquad h=\frac{\dot{r}^2}{2}-\frac{U_0}{r}.
$$
A parabolic solution of the Kepler problem is nothing but a solution with zero energy. There is a unique increasing 
parabolic solution, namely $r(s)=\alpha s^{2/3}$ where $\alpha=(9U_0/2)^{1/3}$. 
Given $0\le a\le b$, the energy of a solution connecting $a$ to $b$ is necessarily greater or 
equal to $-U_0/b$. Moreover, if $0\le a<b$, for $h\ge 0$ or $h=-U_0/b$ there is a unique segment of 
solution of energy $h$ joining $a$ to $b$, this solution increases from $a$ to $b$. 
If $-U_0/b< h<0$ there are exactly two segments of solutions   
of energy $h$ joining $a$ to $b$, a monotonic one, that increases from $a$ to $b$, and
a non-monotonic one, that increases from $a$ to $-U_0/h$ and decreases 
from $-U_0/h$ to $b$.    
Let $\overline{s}(a,b)$ be the time employed by the solution of energy $-U_0/b$ 
to connect $a$ to $b$. We have the following lemma, whose proof is left to the reader.   
\begin{lem} \label{uniq-solution}
Given $0\le a\le b$, and $s>0$, there exists a unique segment of solution joining $a$ to $b$ in time 
$s$, moreover, the solution is monotonic if and only if $0<s\le \overline{s}(a,b)$.  
\end{lem}
\begin{defi}
Given $0\le a\le b$ and $s>0$, we denote by $h(a,b;s)$ the energy of the unique segment of solution joining $a$ 
to $b$ in time $s$, and we denote by $S(a,b;s)$ the Lagrangian action of this solution. 
\end{defi}
\noindent
Since the solution joining $a$ to $b$ in time $s$ is unique, $S(a,b;s)$ is also the minimum of the action
of absolutely continuous paths joining $a$ to $b$ in time $s$.   
\\ \noindent
We shall study the behaviour of the function $r\mapsto h(0,r;s)$ for fixed $s>0$. 
\begin{lem} \label{lem-h-croiss}
Given $s>0$, the function $r\mapsto h(0,r;s)$ is ${\mathcal C}^1$ in $(0,+\infty)$ with a strictly  
positive derivative. Moreover 
$$
\frac{\partial h}{\partial r}(0,\alpha s^\frac{2}{3};s)=\frac{5 U_0}{\alpha^2 s^\frac{4}{3}}.
$$ 
\end{lem}
\noindent
The proof is left to the reader.
We shall also need the following two Propositions   
\begin{prop} \label{S(0,c)}
Let $\overline{\epsilon}>0$. We have 
\begin{equation} \label{ardechoise}
S(0,r;1+\epsilon)= \frac{r^2}{2(1+\epsilon)}+o_r(r^2)
\end{equation}
as $r\rightarrow +\infty$, uniformly for $\epsilon\in [0,\overline{\epsilon}]$. 
\end{prop} 
\Pf\
The parabolic solution $u\mapsto \alpha u^{\frac{2}{3}}$ has zero energy, hence 
$h(0,\alpha(1+\epsilon)^\frac{2}{3};1+\epsilon)=0$. Since we are interested at what 
happens when $r\rightarrow +\infty$, we assume $r> \alpha (1+\overline{\epsilon})^\frac{2}{3}$. 
By Lemma \ref{lem-h-croiss} the energy $h(0,r;1+\epsilon)$ is positive and the solution joining $0$ to $r$
in time $1+\epsilon$ is monotonic.
The function $h=h(0,r;1+\epsilon)$ verifies the identity 
\begin{equation} \label{h-c-eps}
1+\epsilon=\int_0^r \frac{du}{\sqrt{2\left(h+\frac{U_0}{u}\right)}}=\frac{U_0}{2^{\frac{1}{2}} h^{\frac{3}{2}}}
E\left(\frac{hr}{U_0}\right),    
\end{equation} 
where $E:\Rg_+\rightarrow \Rg$ is defined by 
$$
E(x)=\int_0^x \sqrt{\frac{s}{1+s}}ds, \qquad x\in \Rg_+,
$$
and it verifies the estimates
\setlength{\extrarowheight}{0.3cm}
\begin{equation} \label{E-asympt}
\begin{array}{rl}
E(x)&=\frac{2}{3}x^{\frac{3}{2}}+o(x^\frac{3}{2})\quad \text{as}\quad x\rightarrow 0^+ \\ 
E(x)&= x+o(x),\quad \text{as}\quad x\rightarrow +\infty.
\end{array}
\end{equation}
\setlength{\extrarowheight}{0cm}
Let us prove now that 
\begin{equation} \label{h-infini}
h(0,r;1+\epsilon)\rightarrow +\infty,\qquad \text{as}\quad r\rightarrow +\infty
\end{equation}
uniformly on $\epsilon \in [0,\overline{\epsilon}]$. 
Assuming, for the sake of contradiction, that (\ref{h-infini}) is false, there would 
exist two sequence $r_n\rightarrow +\infty$ and $\epsilon_n\in [0,\overline{\epsilon}]$ such that 
$h(0,r_n;1+\epsilon_n)$ is bounded. To simplify notations let us denote 
$h_n=h(0,r_n;1+\epsilon_n)$.
By identities (\ref{h-c-eps}) and (\ref{E-asympt}), the sequence
$h_n r_n$ is bounded too. This implies that $h_n\rightarrow 0$ as
$n\rightarrow +\infty$. Since $E(x)$ is continuous and strictly increasing, 
identity (\ref{h-c-eps}) gives $h_n r_n\rightarrow 0$ as $n\rightarrow +\infty$. 
Applying again (\ref{h-c-eps}) and the first of (\ref{E-asympt}) we obtain 
$$
\lim\limits_{n\rightarrow +\infty} 1+\epsilon_n=\lim\limits_{n\rightarrow +\infty}
\frac{U_0}{\sqrt{2}} \left(\frac{r_n}{U_0}\right)^\frac{3}{2}
\left(\frac{U_0}{h_n r_n}\right)^\frac{3}{2} 
E\left( \frac{h_n r_n}{U_0}\right)=+\infty
$$
that gives a contradiction and proves (\ref{h-infini}).
Writing now (\ref{h-c-eps}) as 
$$
1+\epsilon=\frac{r}{\sqrt{2 h}}\left( \frac{U_0}{hr}\right) 
E\left( \frac{h r}{U_0}\right), 
$$
using the second of (\ref{E-asympt}) we obtain the following estimates 
\begin{equation} \label{h-asympt}
h=h(0,r;1+\epsilon)= \frac{1}{2} \left(\frac{r}{1+\epsilon}\right)^2+o_r(r^2)
\end{equation}
as $r\rightarrow +\infty$, uniformly on $\epsilon\in [0,\overline{\epsilon}]$. 
Let us consider now the action $S(0,r;1+\epsilon)$. 
Let $t\mapsto u(t)$ be the solution joining $0$ with $r$ in time 
$1+\epsilon$. We have 
\setlength{\extrarowheight}{0.4cm}
\begin{equation} \label{S(0,r;1+eps)}
\begin{array}{rl}
S(0,r;1+\epsilon)&=\displaystyle\int_0^{1+\epsilon} \left( \frac{\dot{u}^2}{2}+\frac{U_0}{u}\right) dt=
\displaystyle\int_0^r \frac{h+\frac{2 U_0}{u}}{\sqrt{2\left( h+\frac{U_0}{u}\right)}} du \\
 &=\sqrt{2 h}\displaystyle\int_0^r \sqrt{1+\frac{U_0}{h u}} du - 
\sqrt{\frac{h}{2}}\displaystyle\int_0^r \frac{du}{\sqrt{1+\frac{U_0}{h u}}} \\
&=\frac{U_0}{\sqrt{h}}\left( \sqrt{2} F\left(\frac{hr}{U_0}\right)-\frac{1}{\sqrt{2}} 
E\left(\frac{hr}{U_0}\right) \right),
\end{array}
\end{equation}
\setlength{\extrarowheight}{0cm}
where $F: \Rg_+\rightarrow \Rg$ is defined by  
$$
F(x)=\int_0^x \sqrt{\frac{s+1}{s}} ds,\qquad x\ge 0. 
$$
The function $F$ verifies the asymptotic estimates
\begin{equation} \label{F-asympt}
F(x)= x+o(x),\qquad x\rightarrow +\infty. 
\end{equation}
Replacing (\ref{h-asympt}) in (\ref{S(0,r;1+eps)}) we find (\ref{ardechoise}).
\begin{prop} \label{S(c-point)}
Let $A>0$ and $B>0$ be two constants. 
If we set 
$$
\alpha_0=(8U_0 \alpha)^\frac{1}{2}\quad \text{and}\quad \beta_0=(8U_0)^\frac{1}{2}
$$
then we have 
\begin{equation} \label{barney-panofsky}
S(r,\alpha(s^\frac{2}{3}+\xi);s+\eta)=\alpha_0 s^\frac{2}{3}-\beta_0 r^\frac{1}{2}+o_s(1)
\end{equation}
uniformly on $r\in [0,s^{1/3}]$, $|\xi|\le A$ and $|\eta|\le B$. 
\end{prop}
%
%
%
%
\Pf \
We first prove that the (unique) solution joining $r$ to $\alpha(s^{2/3}+\xi)$ 
in time 
$s+\eta$ is monotonic. In order to simplify the exposition let us term $\lambda(\xi,s)=\alpha(s^{2/3}+\xi)$. We shall 
compare $s+\eta$ with the time employed by
the solution of energy $-\frac{U_0}{\lambda(\xi,s)}$ to connect $r$ to 
$\lambda(\xi,s)$. As usual we denote $\overline{s}(r,\lambda(\xi,s))$ this time.
By definition of $\alpha$  we have 
\setlength{\extrarowheight}{0.5cm}
\begin{displaymath}
\begin{array}{rl}
\overline{s}(r,\lambda(\xi,s))&=\displaystyle\int_r^{\lambda(\xi,s)} 
\frac{du}{\sqrt{2\left(-\frac{U_0}{\lambda(\xi,s)}+\frac{U_0}{u}\right)}} \\
&=\frac{\lambda(\xi,s)^{3/2}}{(2U_0)^{1/2}} 
\displaystyle\int_{\frac{r}{\lambda(\xi,s)}}^1 \frac{dv}{\sqrt{\frac{1}{v}-1}} \\
&=\frac{3s}{2}\left( 1+\frac{\xi}{s^{2/3}}\right)^{3/2} \left( \frac{\pi}{2}-
H\left(\frac{r}{\lambda(\xi,s)}\right) \right),
\end{array}
\end{displaymath}
\setlength{\extrarowheight}{0cm}
where we define  
$$
H :\Rg_+ \rightarrow \Rg,\qquad H(x)=\int_0^x \sqrt{\frac{v}{1-v}} dv.
$$
Since we assume 
$$
0\le r \le s^{1/3} \qquad \text{and}\qquad |\xi|\le A
$$
we have 
$$
\frac{r}{\lambda(\xi,s)}\rightarrow 0\qquad \text{as}\qquad s\rightarrow +\infty.
$$
An easy computation shows that
$$
H(x)=\frac{2}{3}x^{3/2}+{\mathcal O}(x^{5/2}),\qquad x\rightarrow 0,
$$ 
hence we get the estimates 
$$
\overline{s}(r,\lambda(\xi,s))=\frac{3\pi}{4} s \left(1+{\mathcal O}_s(s^{-{1/2}})\right).
$$
Since $\frac{3\pi}{4}>1$, we have $\overline{s}(r,\lambda(\xi,s))>s+\eta$ for $s$ sufficiently
great, and by Lemma \ref{uniq-solution} the solution joining $r$ to $\lambda(\xi,s)$ in time $s+\eta$ 
is monotonic. 
\vspace{1mm} \\
Let $h=h(r,\lambda(\xi,s);s+\eta)$ be the energy of the solution joining $r$ to 
$\lambda(\xi,s)$ in time $s+\eta$. We prove that $h=o_s(1/s)$ for $s\rightarrow +\infty$, 
uniformly on $0\le r\le s^{1/3}$, $|\xi|\le A$ and $|\eta|\le B$. 
The energy $h$ satisfies  the identity 
\setlength{\extrarowheight}{0.5cm}
\begin{equation} \label{rel-h-s-eta-xi}
\begin{array}{rl}
s+\eta&=\displaystyle\int_r^{\lambda(\xi,s)} \frac{ du}{\sqrt{2\left( h+\frac{U_0}{u}\right) }} \\
&=\frac{\lambda(\xi,s)^{3/2}}{(2U_0)^{1/2}}  
\displaystyle\int_{\frac{r}{\lambda(\xi,s)}}^1 
\frac{dv}{\sqrt{\frac{\lambda(\xi,s)}{U_0}h
+\frac{1}{v}}}
\end{array}
\end{equation}
\setlength{\extrarowheight}{0cm}
Introducing the functions
\setlength{\extrarowheight}{0.4cm}
\begin{equation} \label{rocamadour}
\begin{array}{rlrl}
x(r,s,\xi)&=\left(\frac{r}{\lambda(\xi,s)}\right)^{1/2},&\quad
y(s,\xi)&=\frac{\xi}{s^{2/3}}, \\  
k(r,s,\xi,\eta)&=\frac{\lambda(\xi,s)}{U_0} h(r,\lambda(\xi,s);s+\eta),&\quad  
z(s,\eta)&=\frac{\eta}{s},
\end{array}
\end{equation}
\setlength{\extrarowheight}{0cm}
and using the definition of $\alpha$, the relation (\ref{rel-h-s-eta-xi}) becomes
\begin{equation} \label{F(x,y,z,k,c,s,xi,eta)}
F(x(r,s,\xi),y(s,\xi),z(s,\eta),k(r,s,\xi,\eta))=0,
\end{equation}
where $F(x,y,z,k)$ is defined by
$$
F(x,y,z,k)=\int_{x^2}^1 \left(\frac{v}{1+kv}\right)^{1/2}\,dv 
-\frac{2}{3}(1+z)(1+y)^{-{3/2}}.
$$
We think now at $(x,y,z,k)$ as independent variables. 
Using the implicit function theorem we show that the equation 
\begin{equation} \label{F(x,y,z,k)}
F(x,y,z,k)=0
\end{equation}
defines a unique ${\mathcal C}^2$ 
function $k=k(x,y,z)$ for $(x,y,z)$ close to $(0,0,0)$. We observe that $F(x,y,z,k)$ is of class 
${\mathcal C}^2$ with respect to the variables $y$ and $z$. Moreover $F$ is derivable with respect to $x$ and 
$$
\frac{\partial F}{\partial x}(x,y,z,k)=-\frac{2x|x|}{\left(1+kx^2\right)^{1/2}},\qquad \frac{\partial F}{\partial x}(0,0,0,0)=0.
$$ 
$\frac{\partial F}{\partial x}$ is derivable with respect to $x$ and $k$, and we have 
\setlength{\extrarowheight}{0.3cm}
\begin{equation}
\begin{array}{rl}
\frac{\partial^2 F}{\partial x^2}(x,y,z,k)&=-\frac{2|x|(2+kx^2)}{\left(1+kx^2\right)^{3/2}} \\
\frac{\partial^2 F}{\partial k \partial x}(x,y,z,k)&=\frac{x^3 |x|}{\left(1+ kx^2\right)^{3/2}},
\end{array}
\end{equation}
\setlength{\extrarowheight}{0cm}
showing that $\frac{\partial F}{\partial x}$ is of class ${\mathcal C}^1$ in a neighborhood of $(0,0,0,0)$.  
In particular 
$$
\frac{\partial^2 F}{\partial x^2}(0,0,0,0)=0,\qquad \frac{\partial^2 F}{\partial k \partial x}(0,0,0,0)=0.
$$
By the theorem of differentiation under the integral sign, $\frac{\partial F}{\partial k}$, 
$\frac{\partial^2 F}{\partial k^2}$ and $\frac{\partial^2 F}{\partial x \partial k}$ 
are well defined, moreover  
\setlength{\extrarowheight}{0.5cm}
\begin{displaymath}
\begin{array}{rlrl}
\frac{\partial F}{\partial k}(x,y,z,k)&=-\frac{1}{2}\displaystyle\int_{x^2}^1 \left( \frac{v}{1+kv}\right)^{3/2}\, dv,
\qquad 
&\frac{\partial F}{\partial k}(0,0,0,0)&=-\frac{1}{5}, \\
\frac{\partial^2 F}{\partial k^2}(x,y,z,k)&=\frac{3}{4}\displaystyle\int_{x^2}^1 \left(\frac{v}{1+kv}\right)^{5/2}\, dv,
\qquad &\frac{\partial^2 F}{\partial k^2}(0,0,0,0)&=\frac{3}{14}, \\
\frac{\partial^2 F}{\partial x \partial k}(x,y,z,k)&=\frac{x^3 |x|}{\left(1+kx^2\right)^{3/2}}, \qquad
&\frac{\partial^2 F}{\partial x \partial k}(0,0,0,0)&=0.
\end{array}
\end{displaymath}
\setlength{\extrarowheight}{0cm}
By the way, we have also 
$$
\frac{\partial^2 F}{\partial k \partial y}(x,y,z,k)=\frac{\partial^2 F}{\partial k \partial z}(x,y,z,k)=0.
$$
These computations show that $F$ is of class ${\mathcal C}^2$ in a neighborhood of $(0,0,0,0)$. Moreover
$$
F(0,0,0,0)=\int_0^1 \sqrt{v} dv -\frac{2}{3}=0.
$$
By the implicit function theorem, equation (\ref{F(x,y,z,k)}) defines a ${\mathcal C}^2$ function 
$k=g(x,y,z)$ in a neighborhood of $(0,0,0)$ such that $g(0,0,0)=0$ and  
$$
\frac{\partial g}{\partial x}(0,0,0)=\frac{\partial^2 g}{\partial x^2}(0,0,0)=
\frac{\partial^2 g}{\partial x \partial y}(0,0,0)=\frac{\partial^2 g}{\partial x \partial z}(0,0,0)=0,
$$
that is to say
\begin{equation} \label{k-x-y-z}
g(x,y,z)={\mathcal O}(|y|+|z|)+o(x^2+y^2+z^2).
\end{equation}
Coming back to original variables, identity (\ref{k-x-y-z}) gives 
\setlength{\extrarowheight}{0.3cm}
\begin{equation} \label{h-asympt-c-s}
\begin{array}{rl}
h(r,\lambda(\xi,s);s+\eta)&=\frac{U_0}{\lambda(\xi,s)}\, 
g\left( \left(\frac{r}{\lambda(\xi,s)}\right)^{1/2},
\frac{\xi}{s^{2/3}},\frac{\eta}{s}\right) \\
&=o_s({1/s}),
\end{array}
\end{equation}
\setlength{\extrarowheight}{0cm}
as $s\rightarrow +\infty$, uniformly on $0\le r\le s^\frac{1}{3}$, $|\xi|\le A$, and 
$|\eta|\le B$. We compute now the action $S(r,\lambda(\xi,s);s+\eta)$. 
Since the solution joining $r$ to $\lambda(\xi,s)$ in time $s+\eta$ (denoted here $t\mapsto u(t)$) 
is monotonic, we have 
\setlength{\extrarowheight}{0.5cm}
\begin{displaymath}
\begin{array}{rl}
S(r,\lambda(\xi,s);s+\eta)&=\displaystyle\int_0^{s+\eta} 
\left( \frac{\dot{u}^2(t)}{2}+\frac{U_0}{u(t)}\right) dt \\
&=\displaystyle\int_r^{\lambda(\xi,s)} 
\frac{h+\frac{2U_0}{u}}{\sqrt{2\left(h+\frac{U_0}{u}\right)}} du \\
&=\displaystyle\int_r^{\lambda(\xi,s)}\sqrt{2\left(h+\frac{U_0}{u}\right)}\, du-(s+\eta)h.
\end{array}
\end{displaymath}
\setlength{\extrarowheight}{0cm}
Introducing the integration variable $v=\frac{u}{\lambda(\xi,s)}$, by (\ref{h-asympt-c-s}) we find
\begin{equation} \label{sarlat}
S(r,\lambda(\xi,s);s+\eta)=(2 U_0\lambda(\xi,s))^\frac{1}{2}A(x,k)+o_s(1),
\end{equation}
where $x=x(r,s,\xi)$ and $k=k(r,s,\xi,\eta)$ are the functions defined like in (\ref{rocamadour}) and  
$$
A(x,k)=\int_{x^2}^1 \sqrt{k+\frac{1}{v}}\, dv=A_0(k)-B(x,k),
$$
where 
$$
A_0(k)=\int_0^1  \sqrt{k+\frac{1}{v}}\, dv,\qquad B(x,k)=\int_0^{x^2}  \sqrt{k+\frac{1}{v}}\, dv.
$$
Once again, we think at $x$ and $k$ as independent variables and we give an asymptotic expansion of
$A(x,k)$ for $x$ and $k$ close to $0$.
By the classical theorem of differentiation under the integral sign, $A_0(k)$ is derivable in $0$ and 
$$
A_0(k)=2+\frac{k}{3}+o(k).
$$
Moreover we have the following estimates for $B(x,k)$  
\setlength{\extrarowheight}{0.4cm}
\begin{displaymath}
\begin{array}{rl}
B(x,k)&=\displaystyle\int_0^{x^2} \frac{dv}{\sqrt{v}}+\displaystyle\int_0^{x^2} 
\left( \sqrt{k+\frac{1}{v}}-\sqrt{\frac{1}{v}} \right) dv \\
&=2|x|+k\displaystyle\int_0^{x^2} \frac{\sqrt{v}}{\sqrt{1+kv}+1}\, dv \\
&=2|x|+{\mathcal O}(k|x|^3),
\end{array}
\end{displaymath}
\setlength{\extrarowheight}{0cm}
hence
$$
A(x,k)=2+\frac{k}{3}-2|x|+{\mathcal O}(k|x|^3)+o(k),
$$
as $x\rightarrow 0$ and $k\rightarrow 0$.
Replacing in (\ref{sarlat}) and using (\ref{h-asympt-c-s}) we find the final estimates (\ref{barney-panofsky}).
\\ \noindent
The two previous Propositions imply the following one.
\begin{prop} \label{elenuccia}
Given $\overline{\epsilon}>0$, we have
$$
\lim\limits_{\begin{array}{rl} 
             &s\rightarrow +\infty \\ 
             &r\rightarrow +\infty 
             \end{array}} 
{\mathcal N}(r,\alpha s^\frac{2}{3};\epsilon,1,s)=+\infty,
$$ 
uniformly on $\epsilon\in [0,\overline{\epsilon}]$, where ${\mathcal N}$ is the function defined in (\ref{diff-act-kepl}).
\end{prop}
\Pf \
If $0\le \epsilon\le \overline{\epsilon}$ and $0\le r\le s^\frac{1}{3}$, 
from Propositions (\ref{S(0,c)}) and (\ref{S(c-point)}) we have :
$$
{\mathcal N}(r,\alpha s^\frac{2}{3};\epsilon,1,s)=\frac{r^2}{2(1+\epsilon)}(1+o_r(1))-\beta_0 r^\frac{1}{2}+o_s(1),
$$
therefore
\begin{equation} \label{lim-c<s13}
\lim\limits_{\begin{array}{rl} 
             &r\rightarrow +\infty \\ 
             &0\le r \le s^\frac{1}{3}
             \end{array}} 
{\mathcal N}(r,\alpha s^\frac{2}{3};\epsilon,1,s)=+\infty,
\end{equation} 
uniformly on $\epsilon\in [0,\overline{\epsilon}]$. 
Let us consider now the case $r\ge s^\frac{1}{3}$. Forgetting the term
$S(r,\alpha s^\frac{2}{3};s-1)$ in ${\mathcal N}(r,\alpha s^\frac{2}{3};\epsilon,1,s)$ 
and applying again Propositions (\ref{S(0,c)}) and (\ref{S(c-point)}) we find 
$$
\begin{array}{rl}
{\mathcal N}(r,\alpha s^\frac{2}{3};\epsilon,1,s)&\ge \frac{r^2}{2(1+\epsilon)}(1+o_r(1))-\alpha_0 s^\frac{1}{3}+o_s(1) \\
&\ge \frac{s^\frac{2}{3}}{2(1+\epsilon)}(1+o_s(1))-\alpha_0 s^\frac{1}{3}+o_s(1).
\end{array}
$$
This estimates implies the limit 
\begin{equation} \label{lim-c>s13}
\lim\limits_{\begin{array}{rl} 
             &s\rightarrow +\infty \\ 
             &r\ge s^\frac{1}{3}
             \end{array}} 
{\mathcal N}(r,\alpha s^\frac{2}{3};\epsilon,1,s)=+\infty,
\end{equation}
uniformly on $\epsilon\in [0,\overline{\epsilon}]$. 
\\ \noindent
The two limits (\ref{lim-c<s13}) and (\ref{lim-c>s13}) achieve a proof of the Proposition. 
\vspace{3mm}
\begin{center}{\bf Acknowledgements} 
\end{center}
\vspace{4mm}  \noindent
We wish to thanks Alain Chenciner and Albert Fathi for some useful discussion, Marie-Claude Arnaud for  
a careful reading of the manuscript.
%
%
%
%
%
%
%
%
%
%
%
%
\vspace{4mm}\\ 
\noindent
\begin{center}
{\bf References}
\end{center}\par\vspace{5mm}\noindent
[Al] A. Albouy {\it Lectures on the Two-Body Problem, Classical and Celestial Mechanics, 
The Recife Lectures, H. Cabral, F. Diacu editors, Princeton University Press, (2001)}
\par\vspace{1mm} \noindent
[Al-Ch] A. Albouy, A. Chenciner {\it Le probl\` eme des $n$ corps et les
distances mutuelles, Inventiones Mathematic\ae 131, pp. 151-184, (1998)}
\par\vspace{1mm}\noindent
[Ch1] A. Chenciner {\it Collision Totales, Mouvements Compl{\`e}tement
Paraboliques et R{\'e}duction des Homoth{\'e}ties dans le Probl{\`e}me des n Corps, 
Regular and Chaotic Dynamics, Vol. 3, $n^{0}$3, pp. 93-105, (1998)}
\par\vspace{1mm} \noindent 
[Ch2] A. Chenciner {\it Action minimizing solutions of the Newtonian
n-body problem : from homology to symmetry, ICM, Beijing, (2002) } 
\par\vspace{1mm} \noindent
[Cha1] J. Chazy {\it Sur l'allure du mouvement dans le probl\`eme de trois corps quand le temps cro\^it
ind\'efinimment, Ann. Sci. \'Ecole Norm. Sup. 3\`eme s\'erie 39, 29-130, (1922)}
\par\vspace{1mm} \noindent
[Cha2] J. Chazy {\it Sur certaines trajectoires du probl\`eme des $n$ corps, Bulletin Astronomique 35, 321-389, (1918)}
\par\vspace{1mm}\noindent
[Fa] A. Fathi {\it The Weak KAM Theorem in Lagrangian Dynamics, book in preparation}
\par\vspace{1mm} \noindent 
[Fe-Te] D. Ferrario, S. Terracini {\it On the existence of collisionless equivariant 
minimizers for the classical $n$-body problem, Invent. Math.  155, no. 2, 305--362, (2004)}
\par\vspace{1mm}\noindent
[Hu-Sa] N. Hulkower, D. Saari {\it On the manifolds of total collapse orbits and of completely parabolic 
orbits for the $n$-body problem, Journal Diff. Eq. 41, no. 1, 27-43, (1981)}
\par\vspace{1mm}\noindent
[Mad] E. Maderna {\it On weak KAM theory of $N$-body problems, preprint, (2006)}
\par\vspace{1mm} \noindent
[Ma1] C. Marchal {\it How the minimization of action avoids singularities, Celestial Mechanics 83,
pp. 325-354, (2002) }
\par\vspace{1mm} \noindent 
[Ma2] C. Marchal {\it private communication}
\par\vspace{1mm} \noindent
[McG] R. McGehee {\it Triple collisions in the collinear three-body problem, Inv. Math. 27, pp. 191-227, (1974) }
\par\vspace{1mm} \noindent
[Mo] R. Moeckel {\it Orbit Near Triple Collision in the Three-Body Problem, Indiana Univ. Math. Journ., vol. 32, $n^{o}$ 2,
pp. 221-240, (1983)}
\par\vspace{1mm} \noindent
[Ve] A. Venturelli {\it Application de la minimisation de l'action au Probl{\`e}me des 
$N$ corps dans le plan et dans l'espace, Th{\`e}ese de Doctorat, Universit{\'e} de Paris 7-Denis Diderot, 
(2002)} 
\par\vspace{1mm} \noindent

\end{document}